\documentclass[10pt]{article}
\usepackage{amsfonts}
\usepackage{graphicx}
\usepackage{amsmath}
\usepackage{amssymb}
\usepackage{mathtools}
\usepackage{color}
\usepackage{pdflscape}
\usepackage{tabularx}
\usepackage{nicematrix}

\newcommand{\defeq}{\vcentcolon=}

\newtheorem{remark}{Remark}

\makeatother

\begin{document}

\title{Observations on Recurrent Loss in the Neural Network Model of a Partial Differential Equation: the Advection-Diffusion Equation\thanks{This report was prepared as an account of work sponsored by an agency of the United States Government. Neither the United States Government nor any agency thereof, nor any of their employees, make any warranty, express or implied, or assumes any legal liability or responsibility for the accuracy, completeness, or usefulness of any information, apparatus, product, or process disclosed, or represents that its use would not infringe privately owned rights. Reference herein to any specific commercial product, process, or service by trade name, trademark, manufacturer, or otherwise does not necessarily constitute or imply its endorsement, recommendation, or favoring by the United States Government or any agency thereof. The views and opinions of authors expressed herein do not necessarily state or reflect those of the United States Government or any agency thereof.}
}

\author{Jonah A. Reeger %
\thanks{\textit{Email}: jonah.reeger@afit.edu%
} \\
 Department of Mathematics and Statistics\\
 US Air Force Institute of Technology, 2950 Hobson Way\\
 Wright-Patterson AFB, OH 45433, USA\\
}

\maketitle

\begin{abstract}
A growing body of literature has been leveraging techniques of machine learning (ML) to build novel approaches to  approximating the solutions to partial differential equations.  Noticeably absent from the literature is a systematic exploration of the stability of the solutions generated by these ML approaches.  Here, a recurrent network is introduced that matches precisely the evaluation of a multistep method paired with a collocation method for approximating spatial derivatives in the advection diffusion equation.  This allows for two things: 1) the use of traditional tools for analyzing the stability of a numerical method for solving PDEs and 2) bringing to bear efficient techniques of ML for the training of approximations for the action of (spatial) linear operators.  Observations on impacts of varying the large number of parameters in even this simple linear problem are presented.  Further, it is demonstrated that stable solutions can be found even where traditional numerical methods may fail.
\end{abstract}

\section{Introduction}
The solution of time-dependent partial differential equations (PDE) is a topic of concern across many application areas.  Classical approaches to approximating the solution of many PDEs must balance efficiency, accuracy, and stability of the computed solution in the presence of small perturbations to input data and resulting from operations in floating point arithmetic.  Consideration of stability ensures that the behavior of the numerical solution produced by a particular method matches that of an exact solution of a PDE as time increases and without magnifying errors.   These approaches often first construct approximations to any linear differential or integral operators involving space-like variables, which are subject to boundary or some other algebraic conditions.  Such approximations are usually constructed by either expressing the solution as a linear combination of basis functions, e.g., polynomials, trigonometric polynomials, radial kernels, etc., and substituting these expressions directly into the PDE, or by recasting the problem in a weak form and then projecting the dependent variable(s) onto a basis set (see, e.g., \cite{Iserles_2008} for an exposition on the numerical analysis of traditional methods for solving PDEs).  Either way this results in a semi-discrete system of ordinary differential equations (ODEs) that is typically numerically integrated using an appropriate method that considers the stability of the solution.

A growing body of literature has been leveraging techniques of machine learning (ML) to build novel approximations to solutions to PDEs.  Several early works proposed network architectures that either mirrored or incorporated traditional numerical methods in their approach to solve ODEs/PDEs \cite{RICO-MARTÍNEZ01111992,GONZALEZGARCIA1998S965,MEADE19941,712178} or in nonlinear system identification \cite{ALEXANDRIDIS2002479} . More recently, existing methods that pair ML with the solution of PDEs do one of the following:
\begin{itemize}
    \item mimic the method of lines by typically approximating either the action of the differential operators that are present (usually spatial) with a network model or by approximating the entire right-hand side of the PDE through the evaluation of a network (see, e.g., \cite{YBSSHJHMPB2019,math11132791,math12172784}),
    \item utilize a network to encapsulate the entire evolution; that is, evaluation of the neural network itself evolves the solution a single time step, removing the need to involve classical numerical integration (see, e.g., \cite{pmlr-v80-long18a,ZCVCKWDX2022,10.5555/3291125.3291150,RAISSI2019686,LONG2019108925,NORDSTROM2021109821}; also, e.g., \cite{doi:10.1137/24M1686991} in the case of PDEs without time dependence)
\end{itemize}
At this time, the tools that are available for systematically assessing the stability within these ML approaches are limited.  However, it was posited that recurrence in the training of a deep neural network encourage stability in the numerical solution \cite{ZCVCKWDX2022}.  The purpose of this study is to systematically investigate the impacts of recurrence on the construction of weight sets for approximating the action of (spatial) linear differential operators that are used in the process of solving the advection-diffusion equation.

Recurrent Neural Networks (RNNs) can find their origins in the work of McCulloch and Pitts from 1943 \cite{mcculloch43a}, where they considered ``nets with circles" allowing for ``reference to past events of infinite degree of remoteness." Further, Rumelhart, et. al., considered the use of parallel distributed processing in the context of sequential thought processing \cite{Rumelhart86}.  These ideas have an intimate connection to initial value problems, whose solutions over time reference past information.  Several recent works have considered the connection between Neural Networks and dynamic systems \cite{ODERU,LearningDynamical,pmlr-v97-haber19a,Haber} in order to improve the stability, in terms of robustness to perturbations, of the networks for typical tasks (e.g., computer vision, natural language processing, etc.).  In this work the opposite task is investigated: utilize the tools of ML, specifically related to RNNs, to construct stable solutions to a time dependent PDE--the advection-diffusion equation.  A recurrent network is introduced that matches precisely the evaluation of a multistep method paired with a collocation method for approximating spatial derivatives.  This allows for two things: 1) the use of traditional tools for analyzing the stability of a numerical method for solving PDEs and 2) bringing to bear efficient techniques of ML for the training of approximations for the action of (spatial) linear operators.

In order to assist the unfamiliar reader, a comprehensive discussion of the solutions of the advection-diffusion equation, the approximation of the action of linear operators through collocation, the construction of multi-step methods for numerical integration and an introduction to the theory of linear stability of ODEs (and discretized PDEs) is first presented in sections \ref{sec:advecdiff} through \ref{sec:linstab}.  Then, section \ref{sec:networkrep} relates the numerical solution of a linear PDE to a simple neural network, while section \ref{sec:backprop} introduces a recurrent loss function and discusses the procedure for its minimization.  Next, section \ref{sec:numexp} explains the computational experiments that were performed to identify conditions where stable numerical solutions of the advection diffusion equation can be generated through minimization of the loss function.  Observations of the results of these experiments are presented in sections \ref{sec:keyobs} and \ref{sec:otherobs}.  Finally, section \ref{sec:conclusions} concludes this work.  An implementation of the network and training routine discussed in this work is provide on the Matlab File Exchange \cite{ADRNNCode}.

\section{Preliminaries}

This section presents a brief introduction to the theory of linear stability for the unfamiliar reader.  The intent is to introduce an unfamiliar reader to the notion of stability being sought through the numerical procedure in the sections that follow.

\subsection{Asymptotic Behavior of the Advection Diffusion Equation} \label{sec:advecdiff}

Consider the problem of determining the $P$-periodic solution, $u:\mathbb{R}\times[0,T]\rightarrow\mathbb{R}$, of
\begin{subequations}
\label{eq:ContinousProblem}
\begin{align}
\frac{\partial}{\partial t}u(x,t) &= \mathcal{L}u(x,t)\mbox{, }x\in\mathbb{R}\mbox{, }t\in(0,T]\label{eq:ContinuousPDE}\\
u(x,0) &= u_{0}(x)\mbox{, }x\in\Omega,\label{eq:ContinuousIC}
\end{align}
\end{subequations}
where $\mathcal{L}$ is a linear operator acting with respect to $x$ (e.g., a linear combination of various order partial derivatives with respect to $x$). In the case of the advection diffusion equation, which will act as a prototype in this study,
\begin{align}
\mathcal{L} = -c\frac{\partial}{\partial x}+\nu\frac{\partial^{2}}{\partial x^{2}}.\nonumber
\end{align}
This equation has been studied in detail, and comprehensive discussions can be found in many texts on the analysis of PDE (see, e.g., \cite{Guenther_Lee_1996} section 6-5).  Since
\begin{align}
\mathcal{L}\mbox{e}^{i\frac{2\pi \omega}{P}x}=\left(-i\frac{2\pi \omega}{P}c-\left(\frac{2\pi \omega}{P}\right)^{2}\nu\right)\mbox{e}^{i\frac{2\pi \omega}{P}x},\nonumber
\end{align}
with $i=\sqrt{-1}$ reserved here to represent the imaginary unit, a subset of the eigenvalues of $\mathcal{L}$ is (see, e.g., \cite{hunter2001applied})
\begin{align}
\left\{z\in\mathbb{C}:z=-i\frac{2\pi \eta}{P}c-\left(\frac{2\pi \eta}{P}\right)^{2}\nu\mbox{, }\eta\in\mathbb{Z}\right\}.\nonumber
\end{align}
Given that functions that are continuous and $P$-periodic with respect to $x$ can be represented by Fourier series, $u$ can be written as \cite{Guenther_Lee_1996}
\begin{align}
    u(x,t) = \sum\limits_{\eta=-\infty}^{\infty}b_{\eta}(t)\mbox{e}^{i\frac{2\pi \eta}{P}x},\nonumber
\end{align}
which is a linear combination of orthogonal eigenvectors of $\mathcal{L}$.  Further, $P$-periodic solutions to the advection diffusion equation then must satisfy
\begin{align}
    \sum\limits_{\eta=-\infty}^{\infty}\left(\frac{d}{d t}b_{\eta}(t)+\left(i\frac{2\pi \eta}{P}c+\left(\frac{2\pi \eta}{P}\right)^{2}\nu\right)b_{\eta}(t)\right)\mbox{e}^{i\frac{2\pi \eta}{P}x}=0,\nonumber
\end{align}
that is, the time dependent coefficients are defined by
\begin{align}
    b_{\eta}(t) = b_{\eta}(0)\mbox{e}^{-i\frac{2\pi \eta}{P}ct-\left(\frac{2\pi \eta}{P}\right)^{2}\nu t}\nonumber
\end{align}
so that
\begin{align}
    u(x,t) = \sum\limits_{\eta=-\infty}^{\infty}b_{\eta}(0)\mbox{e}^{i\frac{2\pi \eta}{P}(x-ct)-\left(\frac{2\pi \eta}{P}\right)^{2}\nu t}.\label{eq:continuous_solution}
\end{align}
Note that $\{b_{\eta}(0)\}_{n=-\infty}^{\infty}$ is the set of Fourier series coefficients for the initial data.  That is,
\begin{align}
    b_{\eta}(0)=\frac{1}{P}\int\limits_{0}^{P}u_{0}(x)\mbox{e}^{-i\frac{2\pi \eta}{P}x}dx.\label{eq:seriescoefficients}
\end{align}
It is straightforward to show that if $P$-periodic $u$ (with respect to $x$) satisfies the conditions that $\frac{\partial^{p-1}}{\partial x^{p-1}}u(x,t)$ is continuous and $\frac{\partial^{p}}{\partial x^{p}}u(x,t)$ is piecewise continuous for each $x\in\mathbb{R}$, then $b_{\eta}(0)=O(\lvert \eta\rvert^{-p})$ as $\eta\to\pm\infty$ (see, e.g., \cite{Guenther_Lee_1996}).  Further, for such $u$ with $p\geq1$, the convergence implied in \eqref{eq:continuous_solution} is absolute and uniform.  The solution \eqref{eq:continuous_solution} translates the initial data with speed $c$ while the amplitude decays relative to the constant $\nu$.  That is, when the initial data is bounded, solutions to the advection-diffusion equation are at least bounded as $t\to\infty$ and will have amplitude that decays if $\nu>0$.

\subsection{An Approach to Approximating the Action of the Linear Operator} \label{sec:LinOpApprox}

A typical first step in approximating the solution to \eqref{eq:ContinuousPDE} is the consideration of spatial derivatives at each location in a given computational node set $\mathcal{S}=\left\{x_{k}\right\}_{k=0}^{N}\subset\mathbb{R}$.  Since $u$ is periodic with period $P$, it is convenient for $\mathcal{S}$ to be contained entirely in a single period, e.g., $\mathcal{S}\subset[0,P)$.  Further, to reduce computational expense, the approximation at each point is constructed using only a set $\mathcal{N}_{k,n}=\{x_{k,j}\}_{j=1}^{n}$, $n\ll N$, of pairwise distinct points near $x_{k}$.  For reasons that are discussed later, the value of $n$ is chosen to be odd in this work.  It is typical for $\mathcal{N}_{k,n}\subset \mathcal{S}$; however, with the periodicity of $u$ it is useful to artificially extend $\mathcal{S}$ to include two periodic images of the set when selecting information for the local approximations.  That is, two additional sets of nodes, $\mathcal{S}^{\pm}=\left\{x_{k}\pm P\right\}_{k=0}^{N}$, are artificially introduced and $\mathcal{N}_{k,n}$ is then the set of $n$ points in $\mathcal{S}\bigcup \mathcal{S}^{+}\bigcup \mathcal{S}^{-}$ nearest to $x_{k}$.

Now, if $u$ has at least $n$ continuous derivatives with respect to the first argument everywhere in $\Omega$, then for each value of $t$ there exists a polynomial interpolant of degree at most $n-1$ \cite{atkinson}
\begin{align}
    s_{k}(x,t) = \sum\limits_{j=1}^{n}\psi_{k,j}(x)u(x_{k,j},t),\nonumber
\end{align}
that satisfies $s_{k}(x_{k,j},t)=u(x_{k,j},t)$, $j=1,2,\ldots,n$.  The functions $\psi_{k,j}(x)$ are typically Lagrange basis elements constructed on the set $\mathcal{N}_{k,n}$ or some other basis set that satisfies the cardinal property
\begin{align}
    \psi_{k,j}(x_{k,l}) = \left\{\begin{array}{cc} 1, & l=j \\
    0, & l\neq j \end{array}\right..\nonumber
\end{align}
The action of a linear operator, $\mathcal{L}$, acting again with respect to $x$ on $u$ is then approximated locally as
\begin{align}
    (\mathcal{L}u)(x,t) \approx (\mathcal{L}s_{k})(x,t) = \sum\limits_{j=1}^{n}(\mathcal{L}\psi_{k,j})(x)u(x_{k,j},t).\nonumber
\end{align}
Evaluating the local approximation at $x_{k}$ reveals
\begin{align}
    (\mathcal{L}u)(x_{k},t) \approx \sum\limits_{j=1}^{n}w_{k,j}u(x_{k,j},t)\label{eq:localfdapprox}
\end{align}
where $w_{k,j}=(\mathcal{L}\psi_{k,j})(x_{k})$.  In what follows, it will be assumed that the set $\mathcal{S}$ consists of equally spaced points in $[0,P)$.  For instance, a common choice (and the one used herein) is $x_{k} = kh_{x}$, $k=0,1,\ldots,N$, where $h_{x}=P/(N+1)$.  Further, two useful properties are assumed about each set $\mathcal{N}_{k,n}$:
\begin{enumerate}
\item the local computational nodes are arranged in ascending order, i.e., $x_{k,j+1}>x_{k,j}$, $j=1,2,\ldots,n-1$, and
\item the distance, $\lvert x_{k,j}-x_{k}\rvert$, from $x_{k,j}$ to $x_{k}$, $j=1,2,\ldots,n$, is the same for all $k$.
\end{enumerate}
These assumptions on $\mathcal{S}$ and $\mathcal{N}_{k,n}$ ensure that for typical choices of the basis $\{\psi_{k,j}\}_{j=1}^{n}$ the weights satisfy for every $k=1,2,\ldots,N$, $w_{k,j}=w_{j}$, a value independent of $k$.  Specifically, for the advection-diffusion equation the vector of weights defined entrywise as $[\mathbf{w}]_{j}=w_{j}$ can be written as a linear combination of approximations, i.e., $\mathbf{w} = -(c/h_{x}) \mathbf{w}_{1}+(\nu/h_{x}^{2}) \mathbf{w}_{2}$, where the entries of $\mathbf{w}_{1}$ and $\mathbf{w}_{2}$ are the weights for approximating the actions of $\frac{\partial}{\partial x}$ and $\frac{\partial^{2}}{\partial x^{2}}$, respectively, at $\mathbf{x}_{k}$ (scaled to the case of $h_{x}=1$).  To simplify what follows, let
\begin{align}
\boldsymbol{\omega} = \left[\begin{array}{c}\mathbf{w}_{1} \\ \mathbf{w}_{2}\end{array}\right].\nonumber
\end{align}

For a fixed value of $t$, it is typical to approximate the action of $\mathcal{L}$ on $u$ at each point in $\mathcal{S}$ simultaneously through the product of an $N+1\times N+1$ matrix, $D$, with a vector containing values of $u$ evaluated at the nodes in $\mathcal{S}$.  The entries of $D$ are found row by row (something that is easily parallelized).  These entries are defined as
\begin{align}
    [D]_{kl} = \left\{\begin{array}{cc} w_{j} & \mbox{ if }x_{k,j} = x_{l}\mbox{ or }x_{k,j}=x_{l}\pm P\mbox{ for some }(k,j) \\
    0 & \mbox{otherwise}\end{array}\right.,\label{eq:Diffmatentries}
\end{align}
that is, entry $l$ of row $k$ is nonzero only if $x_{l}$ (one of the points in $\mathcal{S}$) or one of its periodic images, $x_{l}\pm P$, appears in the set $\mathcal{N}_{k,n}$ of the points nearest $x_{k}$.  Under the assumptions above placed on $\mathcal{S}$ and $\mathcal{N}_{k,n}$, the matrix $D$ is circulant and sparse and thus has eigenvalues that are easily computed \cite{MatrixComputations}.

Enforcing \eqref{eq:ContinuousPDE} at each point in $\mathcal{S}$ and applying the local approximations for the action of $\mathcal{L}$ leads to the semi-discrete system of equations
\begin{align}
    \frac{d}{d t}\mathbf{u}(t)\approx D\mathbf{u}(t)\label{eq:semidiscrete_vector}
\end{align}
where
\begin{align}
    \mathbf{u}(t) = \left[\begin{array}{cccc}u(x_{0},t) & u(x_{1},t) & \cdots & u(x_{N},t)\end{array}\right]^{T}\nonumber
\end{align}
and
\begin{align}
    \frac{d}{d t}\mathbf{u}(t) = \left[\begin{array}{cccc}\frac{d}{d t}u(x_{0},t) & \frac{d}{d t}u(x_{1},t) & \cdots & \frac{d}{d t}u(x_{N},t)\end{array}\right]^{T}.\nonumber
\end{align}

Further, given that $D$ is a circulant matrix, for any fixed $\psi\in\mathbb{Z}$, it has orthonormal eigenvectors $\mathbf{v}_{\eta+(N+1)\psi}$, $\eta=-\lfloor N/2\rfloor,\ldots,-1,0,1,\ldots,\lfloor N/2\rfloor$, with entries ($k=0,1,2,\ldots,N$ and see, e.g., \cite{MatrixComputations,CirculantMatrices})
\begin{align}
\left[\mathbf{v}_{\eta+(N+1)\psi}\right]_{k+1} = \frac{1}{\sqrt{N+1}}\mbox{e}^{\frac{2\pi i (\eta+(N+1)\psi) k}{N+1}}=\frac{1}{\sqrt{N+1}}\mbox{e}^{\frac{2\pi i \eta k}{N+1}}\nonumber
\end{align}
and corresponding eigenvalues
\begin{align}
    \lambda_{\eta} = \sum\limits_{j=-(n-1)/2}^{(n-1)/2}w_{j+(n+1)/2}\mbox{e}^{\frac{2 \pi i j \eta }{N+1}}.\label{eq:Deigenvalues}
\end{align}
Therefore, when evaluated at $x_{k} = kh_{x}$, $k=0,1,\ldots,N$, (with $h_{x} = P/(N+1)$) continuous, $P$-periodic initial data for \eqref{eq:ContinousProblem} with absolutely convergent series \eqref{eq:continuous_solution} can be expressed in the vector
\begin{align}
    \mathbf{u}(0) =& \sum\limits_{\eta=-\lfloor N/2\rfloor}^{\lfloor N/2\rfloor}\sum\limits_{\psi\in\mathbb{Z}}b_{\eta+(N+1)\psi}(0)\mathbf{v}_{\eta+(N+1)\psi}=\sum\limits_{\eta=-\lfloor N/2\rfloor}^{\lfloor N/2\rfloor}\mathbf{v}_{\eta}\sum\limits_{\psi\in\mathbb{Z}}b_{\eta+(N+1)\psi}(0), \nonumber
\end{align}
with the second equality following from the recognition that $\mathbf{v}_{\eta+(N+1)\psi}=\mathbf{v}_{\eta}$ for all $\psi\in\mathbb{Z}$.  Similarly, the sampled solutions to \eqref{eq:ContinousProblem} can be written
\begin{align}
    \mathbf{u}(t) =\sum\limits_{\eta=-\lfloor N/2\rfloor}^{\lfloor N/2\rfloor}\mathbf{v}_{\eta}\beta_{\eta}(t). \label{eq:EigenvectorSolution}
\end{align}
with
\begin{align}
\beta_{\eta}(t) = \sum\limits_{\psi\in\mathbb{Z}}b_{\eta+(N+1)\psi}(0)\mbox{e}^{-i\frac{2\pi (\eta+(N+1)\psi)}{P}ct-\left(\frac{2\pi (\eta+(N+1)\psi)}{P}\right)^{2}\nu t} \nonumber
\end{align}
That is, it is likely (dependent on $b_{\eta+(N+1)\psi}(0)$) that some of $\mathbf{u}(t)$ is in the direction of each of the eigenvectors of $D$.

\subsection{Time Integration with Multistep Methods}

The semidiscrete PDE can be numerically integrated with, for instance, a multistep method as (see, e.g., \cite{atkinson})
\begin{align}
    \tilde{\mathbf{u}}(t+h_{t}) = \tilde{\mathbf{u}}(t)+h_{t}\sum_{l=0}^{s-1}\alpha_{l}D\tilde{\mathbf{u}}(t-lh_{t}), \label{eq:ms_method}
\end{align}
with $h_{t}>0$ and the tilde indicating the values of the numerical solution.  To simplify the ensuing discussion consider, for now, the scalar ordinary differential equation (ODE)
\begin{subequations}\label{eq:simple_ODE}
\begin{align}
    \frac{d}{dt}\mu(t) =& f(t,\mu(t)), \mbox{ }t>0\\
     \mu(0) =& \mu_{0}
\end{align}
\end{subequations}
the set of coefficients $\{\alpha_{l}\}_{l=0}^{s-1}$ in \eqref{eq:ms_method} is constructed by way of
\begin{align}
    \mu(t+h_{t}) &= \mu(t)+\int\limits_{t}^{t+h_{t}}f(\tau,\mu(\tau))d\tau\nonumber\\
    &\approx \mu(t)+h_{t}\sum_{l=0}^{s-1}\alpha_{l} f(t-lh_{t},\mu(t-lh_{t})),\label{eq:ms_integral}
\end{align}
where
\begin{align}
\alpha_{l} = \int\limits_{0}^{1}\prod\limits_{\begin{array}{c}j=0\\ j\neq(s-l)\end{array}}^{s-1}{\frac{\tau+s-(j+1)}{s-l-j}}d\tau\label{eq:ms_Weights}
\end{align}
results from integrating a Lagrange basis polynomial now used in the interpolation of $f$ at the equally spaced values of $t-lh_{t}$, $l=0,1,\ldots,s-1$.  Consequently, these methods, known as $s$-step Adams-Bashforth methods (AB-$s$), are explicit and they are exact when $f$ is a polynomial of degree at most $s-1$ with respect to its first argument.  For more generic and sufficiently smooth $v$, when exact values of $\mu(t-lh_{t})$, $l=0,1,\ldots,s-1$, are known this approximation is $O(h_{t}^{s+1})$ accurate (in terms of truncation error) for $\mu(t+h_{t})$, as $h_{t}\rightarrow0$.  Considering approximate numerical integration in time, this choice leads to a method given by
\begin{align}
    \tilde{\mu}(t+h_{t}) = \tilde{\mu}(t)+h_{t}\sum_{l=0}^{s-1}\alpha_{l} f(t-lh_{t},\tilde{\mu}(t-lh_{t})) \label{eq:one_d_ms_methd}
\end{align}
that is $O(h_{t}^{s})$ accurate \cite{atkinson}.

\subsection{Linear Stability} \label{sec:linstab}

In the further simplified linear case of \eqref{eq:simple_ODE}, where $f(t,\mu(t))=z \mu(t)$ for some scalar $z\in\mathbb{C}$, it is common to analyze for which values of the quantity $h_{t}z$ the numerical solution of \eqref{eq:simple_ODE} generated by the method \eqref{eq:one_d_ms_methd} match the asymptotic behavior of the exact solution $\mu(t)$ as $t\to\infty$ (or at least remain bounded).  When $z$ has nonpositive real part, the exact solution to \eqref{eq:simple_ODE}, given by $\mu(t) = \mu_{0}\exp(z t)$, is bounded and when the real part of $z$ is negative, $\mu(t)\to0$ as $t\to \infty$.  In the case of $z$ with positive real part, solutions grow exponentially, which is a difficulty for many numerical methods.

Notice that
\begin{align}
    \tilde{\mu}(t+h_{t}) = \tilde{\mu}(t)+h_{t}\sum_{l=0}^{s-1}\alpha_{l} z\tilde{\mu}(t-lh_{t}) \label{eq:one_d_ms_methd}
\end{align}
is a linear recurrence relation.  It is helpful to note that this recurrence relation appears as the first equation of the system of linear equations
\begin{align}
    \tilde{\boldsymbol{\mu}}(t+h_{t}) = M \tilde{\boldsymbol{\mu}}(t).\nonumber
\end{align}
where $\tilde{\boldsymbol{\mu}}(t) = \left[\begin{array}{cccc}\tilde{\mu}(t) & \tilde{\mu}(t-h_{t}) & \cdots & \tilde{\mu}(t-(s-1)h_{t})\end{array}\right]^{T}$ and
\begin{align}
    M = \begin{bNiceArray}{cccccc}
    1+\alpha_{0}h_{t}z & \alpha_{1}h_{t}z & \alpha_{2}h_{t}z & \cdots & \alpha_{s-2}h_{t}z & \alpha_{s-1}h_{t}z \\
    1 & 0 & 0 & \cdots & 0 & 0\\
    0 & 1 & 0 & \cdots & 0 & 0 \\
     & 0 & \ddots & \ddots & \vdots & \vdots \\
    & & \ddots & 1 & 0 & 0\\
      & & & 0 & 1 & 0
\end{bNiceArray}.\nonumber
\end{align}

As long as $\tilde{\mu}(jh_{t})$, $j=0,1,\ldots,s-1$, are known values it is straightforward to reason that for $l>0$
\begin{align}
    \tilde{\boldsymbol{\mu}}((s+l)h_{t}) = M^{l+1} \tilde{\boldsymbol{\mu}}((s-1)h_{t}).\nonumber
\end{align}
The matrix $M$ has characteristic polynomial
\begin{align}
p(\zeta) = (-1)^{s}(\zeta^{s}-(1+\alpha_{0}h_{t}z)\zeta^{s-1}-\sum\limits_{j=1}^{s-1}\alpha_{j}h_{t}z \zeta^{s-(j+1)}),\nonumber
\end{align}
whose roots, $\zeta_{0},\zeta_{1},\ldots,\zeta_{s}$, are the eigenvalues of $M$. Further,
\begin{align}
\boldsymbol{\zeta}_{j} = \left[\begin{array}{ccccc}\zeta_{j}^{s-1} & \zeta_{j}^{s-2} & \cdots & \zeta_{j} & 1\end{array}\right]^{T}\nonumber
\end{align}
is an eigenvector of $M$ associated with the eigenvalue $\zeta_{j}$. 
Let $M = W(Z+P)W^{-1}$ be the Jordan form of $M$ where $Z$ is the diagonal matrix with $[Z]_{jj}=\zeta_{j}$ and $P$ is a nilpotent matrix with $\lVert P\rVert=1$.  Then
\begin{align}
    \tilde{\boldsymbol{\mu}}((s+l)h_{t}) = W(Z+P)^{l+1}W^{-1} \tilde{\boldsymbol{\mu}}((s-1)h_{t})\mbox{, }l>0.\nonumber
\end{align}
From this it is clear that $\lVert \tilde{\boldsymbol{\mu}}((s+l)h_{t})\rVert$ is bounded as $l$ increases if $|\zeta_{i}|\leq1$ when $\zeta_{i}$ is a simple root of $p$ or is not defective as an eigenvalue of $M$ and if $|\zeta_{i}|<1$ for a root of greater multiplicity that is defective \cite{Leveque}.  Expressing $\zeta$ in polar form, evaluation of (as in, e.g., \cite{MGBFJAR2014})
\begin{align}
    z h = \frac{\mbox{e}^{is\theta}-\mbox{e}^{i(s-1)\theta}}{\sum_{l=0}^{s-1}\alpha_{l}\mbox{e}^{-i(s-1-l)\theta}}\nonumber
\end{align}
with $\theta\in[0,2\pi)$ traces a curve in the complex plane, indicating the choices of $\xi = z h_{t}$ for which $|\zeta|=1$.  The linear stability region corresponds to those where $\lVert \tilde{\boldsymbol{\mu}}((s+l)h_{t})\rVert$ is bounded.

Returning to the solving \eqref{eq:semidiscrete_vector}, since $D$ is unitarily diagonalizable, \eqref{eq:ms_method} can be expressed as
\begin{align}
    V^{*}\tilde{\mathbf{u}}(t+h_{t}) = V^{*}\tilde{\mathbf{u}}(t)+h_{t} \sum_{l=0}^{s-1}\alpha_{l}\Lambda V^{*}\tilde{\mathbf{u}}(t-lh_{t}), \nonumber
\end{align}
with $V=\left[\begin{array}{cccc}\mathbf{v}_{-\lfloor N/2\rfloor},\ldots,\mathbf{v}_{-1},\mathbf{v}_{0},\mathbf{v}_{1},\ldots,\mathbf{v}_{\lfloor N/2\rfloor}\end{array}\right]$ and $\Lambda$ the diagonal matrix with $[\Lambda]_{kk}=\lambda_{k-\lfloor N/2\rfloor-1}$.  Therefore, each component of $V^{*}\tilde{\mathbf{u}}(t)$ should have the asymptotic behavior of the corresponding component of $V^{*}\mathbf{u}(t)$.  Further, it is clear from \eqref{eq:EigenvectorSolution} that for linear stability (considering that solutions to \eqref{eq:ContinousProblem} are bounded for the advection diffusion equation) values of $h_{t}$ such that $h_{t}\lambda$ is contained in the linear stability region of the method should be sought for each eigenvalue, $\lambda_{\eta}$, of $D$.

\section{A Network Representation of a Single Step} \label{sec:networkrep}

\begin{figure}
\begin{center}
\includegraphics[width=\linewidth]{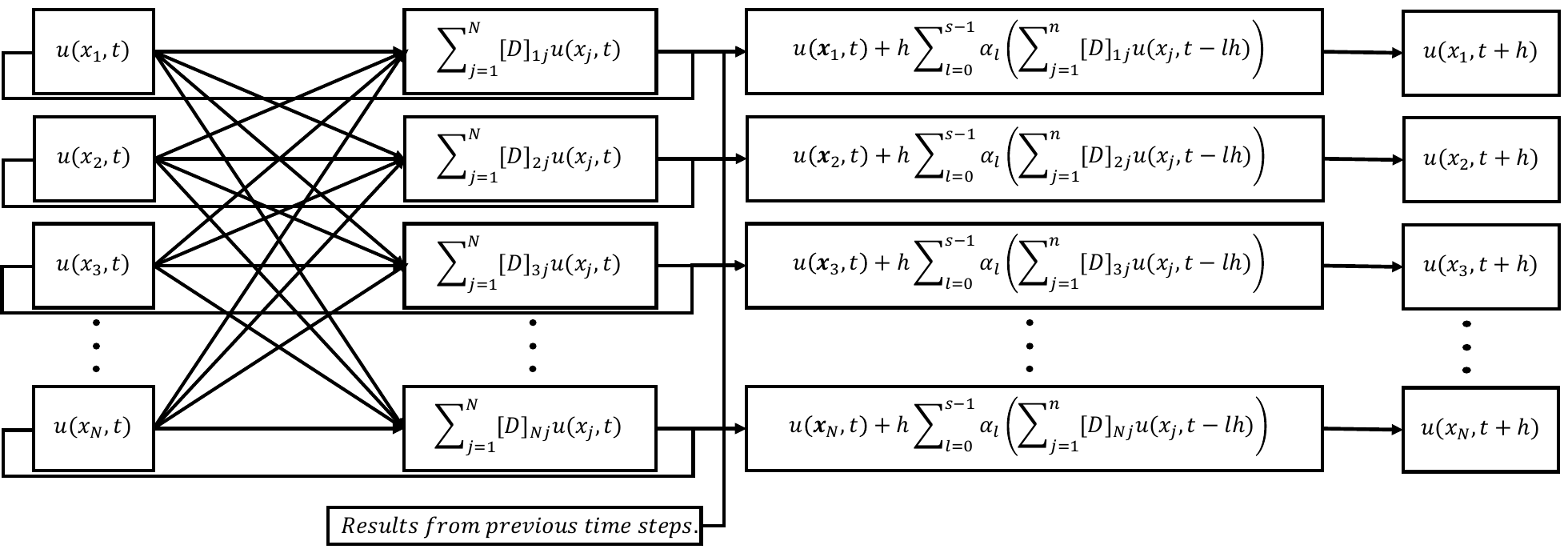}
\end{center}
\caption{A simplified schematic of a single step of an $s$-step multistep method for a semi-discrete PDE cast as a network.}
\label{fig:network_schematic}
\end{figure}

Departing from traditional approaches to the development of methods for solving \eqref{eq:ContinousProblem}, the step from $t$ to $t+h_{t}$ via \eqref{eq:ms_method} can alternatively be viewed as the evaluation of a neural network, albeit one with linear activation functions and no bias.  A simple schematic of this single step can be seen in figure \ref{fig:network_schematic}.  Under the assumptions that the computational node set is equally spaced and the nodes in $\mathcal{N}_{k,n}$, $k=0,1,\ldots,N$, are sorted in ascending order, here the network is considered to, instead, have five separate layers to assist in the use of backpropagation to compute the gradients of the cost function defined below.

For $l\geq 0$, the layers are evaluated, in turn, beginning with the input layer.  The output of each layer is collected in an intermediate variable denoted by either $\boldsymbol{z}$ or $Z$, depending on whether the output is a vector or matrix, respectively.  These intermediate variables have a superscript indicating which layer they are the output of.   When considering training the network, the subscript $\tau$ is convenient in indicating that these outputs change with the training data (which are enumerated by $\tau$).  Otherwise, when evaluating the network outside of the training context the superscript could be ignored or dropped.

The input layer (superscript $I$) has $N+1$ nodes which can be represented by the $(N+1)\times 1$ vector
\begin{align}
{\boldsymbol{z}}_{\tau}^{I}((s-1+l)h_{t};\boldsymbol{\omega}) =\left\{\begin{array}{cc}\tilde{\mathbf{u}}((s-1+l)h_{t};\boldsymbol{\omega}), & l = 0 \\
{\boldsymbol{z}}_{\tau}^{+}((s-1+l-1)h_{t};\boldsymbol{\omega}), & l>0\end{array}\right..\nonumber
\end{align}
Here and in what follows, the dependence of these values on $\boldsymbol{\omega}$, defined in section \ref{sec:LinOpApprox}, is expressed clearly since the training of the neural network will seek the weights represented by this vector specifically when minimizing the objective described in the next section.  While the entries of the vector $\boldsymbol{\omega}$ are used as decision variables in this work, it is also possible to include the coefficients of the multi-step method, $\alpha_{l}$, $l=0,1,\ldots,s-1$, as parameters that can be optimized when attempting to construct a numerical method that is stable.  This approach has been attempted by the author, with some limited success, but has been omitted from this work for brevity and since fixing the weights of the multi-step method produced more consistent results.

Notice that when $l=0$, the input layer is simply the value of $\tilde{\mathbf{u}}((s-1)h_{t};\boldsymbol{\omega})$.  This and the values of $\tilde{\mathbf{u}}(kh_{t};\boldsymbol{\omega})$, $k=0,1,\ldots,s-2$ are assumed to be available to kickstart the method, which is a common requirement for multistep methods.  These values could be obtained using a one-step method of appropriate order, such as a Runge-Kutta or Taylor series method.  For $l>0$, the values ${\boldsymbol{z}}_{\tau}^{I}((s-1+l)h_{t};\boldsymbol{\omega})$, and so ${\boldsymbol{z}}_{\tau}^{+}((s-1+l-1)h_{t};\boldsymbol{\omega})$, are simply values of $\tilde{\mathbf{u}}((s-1+l)h_{t};\boldsymbol{\omega})$ computed by a feedforward pass through the network.

Next, the assignment layer (superscript $A$) duplicates and arranges the values in ${\boldsymbol{z}}_{\tau}^{I}((s-1+l)h_{t})$ into the $n(N+1)\times 1$ vector
\begin{align}
    {\boldsymbol{z}}_{\tau}^{A}((s-1+l)h_{t};\boldsymbol{\omega}) = U{\boldsymbol{z}}^{I}((s-1+l)h_{t};\boldsymbol{\omega}),\nonumber
\end{align}
with $U$ a constant $n(N+1)\times N+1$ sparse matrix defined by
\begin{align}
    [U]_{ij}=\left\{\begin{array}{cc}1, & x_{j}\in\mathcal{N}_{k,n}\mbox{ and }n(k-1)+1\leq i\leq nk \\ 0, & \mbox{otherwise}\end{array}\right..\nonumber
\end{align}
The first $n$ entries of ${\boldsymbol{z}}_{\tau}^{A}((s-1+l)h_{t};\boldsymbol{\omega})$ are the elements of ${\boldsymbol{z}}_{\tau}^{I}((s-1+l)h_{t};\boldsymbol{\omega})$ that correspond to $\mathcal{N}_{1}$ (in ascending order), the next $n$ entries are the same for $\mathcal{N}_{2}$, and so on.  That is,

\begin{align}
    {\boldsymbol{z}}_{\tau}^{A}((s-1+l)h_{t};\boldsymbol{\omega}) = \begin{bNiceArray}{c}\tilde{u}(x_{0,1},(s-1+l)h_{t};\boldsymbol{\omega})\\\tilde{u}(x_{0,2},(s-1+l)h_{t};\boldsymbol{\omega})\\\vdots\\\tilde{u}(x_{0,n},(s-1+l)h_{t};\boldsymbol{\omega})\\\hdottedline
    \tilde{u}(x_{1,1},(s-1+l)h_{t};\boldsymbol{\omega})\\\tilde{u}(x_{1,2},(s-1+l)h_{t};\boldsymbol{\omega})\\\vdots\\\tilde{u}(x_{1,n},(s-1+l)h_{t};\boldsymbol{\omega})\\\hdottedline\vdots\\\hdottedline
    \tilde{u}(x_{N,1},(s-1+l)h_{t};\boldsymbol{\omega})\\\tilde{u}(x_{N,2},(s-1+l)h_{t};\boldsymbol{\omega})\\\vdots\\\tilde{u}(x_{N,n},(s-1+l)h_{t};\boldsymbol{\omega})\end{bNiceArray}.\nonumber
\end{align}

Third, the reshape layer (superscript $R$) constructs a full $(N+1)\times n$ matrix
\begin{align}
Z_{\tau}^{R}((s-1+l)h_{t};\boldsymbol{\omega}) = \sum\limits_{j=1}^{n}R_{j}{\boldsymbol{z}}_{\tau}^{A}((s-1+l)h_{t};\boldsymbol{\omega})\mathbf{e}_{j}^{T}\nonumber
\end{align}
from the values of ${\boldsymbol{z}}_{\tau}^{A}((s-1+l)h_{t};\boldsymbol{\omega})$, with $R_{j}$, $j=1,2,\ldots,n$, a set of constant $(N+1)\times n(N+1)$ sparse matrices containing entries
\begin{align}
    [R_{j}]_{il} = \left\{\begin{array}{cc} 1 & l=(i-1)n+j
    \\ 0 & \mbox{otherwise}\end{array}\right..\nonumber
\end{align}
Note that $\left[Z_{\tau}^{R}((s-1+l)h_{t};\boldsymbol{\omega})\right]_{kj}$ corresponds to $x_{k,j}\in\mathcal{N}_{k,n}$, i.e.,
\begin{align}
    Z_{\tau}^{R}(\cdot;\boldsymbol{\omega}) = \left[\begin{array}{cccc}\tilde{u}(x_{0,1},\cdot;\boldsymbol{\omega}) & \tilde{u}(x_{0,2},\cdot;\boldsymbol{\omega}) & \cdots & \tilde{u}(x_{0,n},\cdot;\boldsymbol{\omega})\\
    \tilde{u}(x_{1,1},\cdot;\boldsymbol{\omega}) & \tilde{u}(x_{1,2},\cdot;\boldsymbol{\omega}) & \cdots & \tilde{u}(x_{1,n},\cdot;\boldsymbol{\omega})\\
    \vdots & \vdots & & \vdots\\
    \tilde{u}(x_{N,1},\cdot;\boldsymbol{\omega}) & \tilde{u}(x_{N,2},\cdot;\boldsymbol{\omega}) & \cdots & \tilde{u}(x_{N,n},\cdot;\boldsymbol{\omega})\end{array}\right].\nonumber
\end{align}

Next, the convolution layer (superscript $C$) performs the operation $D\tilde{\mathbf{u}}((s-1+1)h_{t};\boldsymbol{\omega})$ through the matrix/vector multiplication defined in
\begin{align}
{\boldsymbol{z}}_{\tau}^{C}((s-1+l)h_{t};\boldsymbol{\omega}) = Z_{\tau}^{R}((s-1+l)h_{t};\boldsymbol{\omega})\left(-\frac{c}{h_{x}}\mathbf{w}_{1}+\frac{\nu}{h_{x}^{2}}\mathbf{w}_{2}\right),\nonumber
\end{align}
so that ${\boldsymbol{z}}_{\tau}^{C}(\cdot;\boldsymbol{\omega})$ is an $N+1\times 1$ vector when both $\mathbf{w}_{1}$ and $\mathbf{w}_{2}$ are $n\times 1$ vectors.

Finally, the additive layer (superscript $+$) evaluates the multistep method \eqref{eq:ms_method} by way of
\begin{align}
    {\boldsymbol{z}}_{\tau}^{+}((s-1+l)h_{t};\boldsymbol{\omega}) = {\boldsymbol{z}}_{\tau}^{I}((s-1+l)h_{t};\boldsymbol{\omega})+h_{t}\sum\limits_{j=0}^{s-1}\alpha_{j}{\boldsymbol{z}}_{\tau}^{C}((s-1+l-j)h_{t};\boldsymbol{\omega})\nonumber
\end{align}
Note that for $l=0,1,\ldots,s-1$, this expression relies on the values of ${\boldsymbol{z}}_{\tau}^{C}(kh_{t};\boldsymbol{\omega})$, $k=0,1,\ldots,s-1$, which have not been computed.  These correspond to the evaluation of $D\tilde{\mathbf{u}}(kh_{t};\boldsymbol{\omega})$, which can be completed directly by multiplication or by evaluating in sequence ${\boldsymbol{z}}_{\tau}^{I}(kh_{t};\boldsymbol{\omega}) = \tilde{\mathbf{u}}(kh_{t};\boldsymbol{\omega})$, ${\boldsymbol{z}}_{\tau}^{A}(kh_{t};\boldsymbol{\omega}) = U{\boldsymbol{z}}^{I}(kh_{t};\boldsymbol{\omega})$, $Z_{\tau}^{R}(kh_{t};\boldsymbol{\omega}) = \sum\limits_{j=1}^{n}R_{j}{\boldsymbol{z}}_{\tau}^{A}(kh_{t};\boldsymbol{\omega})\mathbf{e}_{j}^{T}$ and ${\boldsymbol{z}}_{\tau}^{C}(kh_{t};\boldsymbol{\omega}) = Z_{\tau}^{R}(kh_{t};\boldsymbol{\omega})\left(-\frac{c}{h_{x}}\boldsymbol{\omega}_{1}+\frac{\nu}{h_{x}^{2}}\boldsymbol{\omega}_{2}\right)$.  The latter is preferred because it more easily highlights the dependence of ${\boldsymbol{z}}_{\tau}^{C}(kh_{t};\boldsymbol{\omega})$ on the entries of $\boldsymbol{\omega}$ for backpropagation.  A time-stepping method traverses this network many times, with results from each step clearly impacting all later times, and creating a nonlinear dependence on the parameter set, $\boldsymbol{\omega}$.

\begin{remark}
While the network architecture presented here includes only linear activation functions, so that a feedforward pass through the network is analogous to the evaluation of a multi-step method for the linear PDE considered here, nonlinear activation functions could be included without much change to the architecture (or to the backpropagation procedure discussed next).  Further, this architecture could be extended to the case of a nonlinear PDE, even with only linear activation functions up to the convolution layer.  This could appear as an evaluation of the right hand side of a nonlinear PDE after the convolution layer approximates the action of any linear operators that are present utilizing the trained weights.  The evaluation of the nonlinear right hand side would be comparable to the evaluation of a nonlinear activation function.  Likewise, if a time integrator like a Runge-Kutta (RK) method were used in place of the multi-step method, then the architecture could be adjusted to include multiple recursive passes from the input layer to the convolution layer to populate the required stages before evaluation of an appropriate additive layer.  The inputs to these recursive passes would be combinations of the outputs of the previous passes along with scalar multiples of $h_{t}$, all of which would be defined by the particular RK method in use.
\end{remark}


\section{Recurrent Loss and Backpropagation} \label{sec:backprop}
The cost (loss) function considered here is (similar to that of \cite{ZCVCKWDX2022})
\begin{align}
J(\boldsymbol{\omega}) = \sum\limits_{\tau=1}^{T}\sum\limits_{q=0}^{Q-1}J_{\tau q}(\boldsymbol{\omega}),\nonumber
\end{align}
with
\begin{align}
    J_{\tau q}(\boldsymbol{\omega})=\frac{1}{2}\left\lVert {\boldsymbol{z}}_{\tau}^{+}((s-1+q)h_{t};\boldsymbol{\omega})-\mathbf{u}_{\tau}((s+q)h_{t})\right\rVert_{2}^{2}.\nonumber
\end{align}
In this expression for the cost, $\tau$ enumerates $T$ separate cases of initial conditions that are used to generate training data, $\mathbf{u}_{\tau}(t)$, and $q$ counts $Q$ time steps beyond the initial $s$ steps required to kickstart the multi-step method beginning at $t=0$.  The intention is to choose $Q>1$ to promote stability in the solution of the PDE by minimizing the forward error with respect to entries of $\boldsymbol{\omega}$ across multiple steps at once.  Due to the small number of decision variables here, a quasi-Newton method (which requires computation of the gradient and Hessian of $J$ with respect to these variables) is favored in determining the extreme points of $J$.  This choice of an iterative method is further justified given that ${\boldsymbol{z}}_{\tau}^{+}((s-1+q)h_{t};\boldsymbol{\omega})$ (certainly for $Q>1$) depends exponentially on the decision variables, since equation \eqref{eq:ms_method} identifies this vector relying on the repeated application of $D$.   Given a particular value of $q$, call it $q^{*}$, the value of ${\boldsymbol{z}}_{\tau}^{+}((s-1+q^{*})h_{t};\boldsymbol{\omega})$ is dependent on all values of ${\boldsymbol{z}}_{\tau}^{+}((s-1+q)h_{t};\boldsymbol{\omega})$ for which $q<q^{*}$ and so all must be considered, in turn, when computing the gradient (and hessian) of the loss.  The results for each value of $\tau$ are independent from one another, so when minimizing $J$ it suffices to consider each $J_{\tau q}$ in turn.

As is typical for the computation of the gradient of the objective via backpropagation, it is useful to compute a variety of intermediate quantities that leverage the chain rule of differentiation.  These quantities are summarized in tables \ref{tab:backprop1} and \ref{tab:backprop2} and they are represented with the variables $\boldsymbol{\delta}$ or $\Delta$, depending on whether the variable is a vector or matrix, respectively.  The superscript again refers to the layer of the network for which the intermediate quantity is computed, while the subscripting now indicates the case of training data ($\tau$) and the time step in the recurrent loss ($q$).  These quantities can be used to compute the gradient with respect to the weights as
\begin{align}
\nabla_{\boldsymbol{\omega}}J(\boldsymbol{\omega}) = \left[\begin{array}{c}-\frac{c}{h_{x}}I \\ \frac{\nu}{h_{x}^{2}}I\end{array}\right]\sum\limits_{\tau=1}^{T}\sum\limits_{q=0}^{Q-1}\sum\limits_{l=0}^{s-1+q} Z_{\tau}^{R}(lh_{t};\boldsymbol{\omega}){\boldsymbol\delta}_{\tau q}^{C}(lh_{t};\boldsymbol{\omega}).\nonumber
\end{align}
Here and in tables \ref{tab:backprop1} and \ref{tab:backprop2}, $\nabla$ represents the gradient operator with the subscript indicating the variables for which the partial derivatives should be computed.  For instance, $\nabla_{\boldsymbol{\omega}}=\left[\begin{array}{cccc}\partial/\partial[\boldsymbol{\omega}]_{1} & \partial/\partial[\boldsymbol{\omega}]_{2} & \cdots & \partial/\partial[\boldsymbol{\omega}]_{2n}\end{array}\right]$.

The quasi-Newton method attempts to find roots of the gradient of $J$ with respect to the decision variables.  Beginning with an initial guess at the weights, $\boldsymbol{\omega}^{(0)}$, the approximations for the roots are updated as
\begin{align}
    \boldsymbol{\omega}^{(\kappa+1)} = \boldsymbol{\omega}^{(\kappa)}+\rho^{(\kappa)}\boldsymbol{\sigma}^{(\kappa)},\label{eq:line_search_step}
\end{align}
with $\rho^{(\kappa)}\in(0,1]$ a line search parameter and $\boldsymbol{\sigma}^{(\kappa)}$ satisfying
\begin{align}
    H^{(\kappa)}\boldsymbol{\sigma}^{(\kappa)}=-\nabla_{\boldsymbol{\omega}}J^{(\kappa)}(\boldsymbol{\omega})\big\rvert_{\boldsymbol{\omega}=\boldsymbol{\omega}^{(\kappa)}}. \nonumber
\end{align}
The superscript in parentheses, with $\kappa = 0,1,2,\ldots$, in these expressions indicates the iteration, or epoch in the language of neural networks, of the method for minimizing the cost. Therefore, $H^{(\kappa)}$ is meant to represent the Hessian of $J$ (or an approximation thereof) with respect to $\boldsymbol{\omega}$ and evaluated at $\boldsymbol{\omega}^{(\kappa)}$.   Instead of computing $H^{(\kappa)}$ exactly at each iteration, which is significantly more expensive than backpropagation for the gradient, the BFGS approximation \cite{Nocedal_Wright_2006} is utilized instead.  The approximate Hessian is initialized as
\begin{align}
H^{(0)}=\left\lVert\nabla_{\boldsymbol{\omega}}J(\boldsymbol{\omega})\big\rvert_{\boldsymbol{\omega}=\boldsymbol{\omega}^{(0)}}\right\rVert_{2}I\nonumber
\end{align}
$I$ is the identity matrix of appropriate size.

Note that the line search parameter $\rho^{(\kappa)}$ \cite{Nocedal_Wright_2006} chosen so that
\begin{align}
J(\boldsymbol{\omega}^{(\kappa)}+\rho^{(\kappa)}\boldsymbol{\sigma}^{(\kappa)})<J(\boldsymbol{\omega}^{(\kappa)}).\label{eq:line_search_condition}
\end{align}
At iteration $\kappa$, $\rho^{(\kappa)}$ is initialized as one, and $J$ is evaluated at the result of \eqref{eq:line_search_step}.   While $J(\boldsymbol{\omega}^{(\kappa)}+\rho^{(\kappa)}\boldsymbol{\sigma}^{(\kappa)})\geq J(\boldsymbol{\omega}^{(\kappa)})$, $\rho^{(\kappa)}$ is successively decreased by a multiplicative factor of $3/4$.  This process continues until \eqref{eq:line_search_condition} is met, at which point $\boldsymbol{\omega}^{(\kappa+1)}$ is assigned by \eqref{eq:line_search_step}.  If $\rho^{(\kappa)}$ decreases below a chosen lower bound (here $10^{-5}$), then the optimization algorithm is terminated and assumed to have failed.

\begin{landscape}
\begin{table}[p]
\begin{center}
\begin{tabularx}{1.05\linewidth}{|c|X|}
\hline
\begin{tabular}{c}Additive \\ Layer\end{tabular} &
\begin{minipage}{1.05\linewidth}
\begin{equation}
{\boldsymbol\delta}_{\tau q}^{+}((s-1+l)h_{t}) \defeq \nabla_{{\boldsymbol{z}}_{\tau}^{+}((s-1+l)h_{t};\boldsymbol{\omega})} J_{\tau q}(\boldsymbol{\omega}) = \left\{\begin{array}{cc}{\boldsymbol{z}}_{\tau}^{+}((s-1+q)h_{t};\boldsymbol{\omega})-\mathbf{u}_{\tau}((s+q)h_{t}), & l=q \\
{\boldsymbol\delta}_{\tau q}^{I}((s+l)h_{t};\boldsymbol{\omega}), & 0\leq l< q\end{array}\right.\nonumber
\end{equation}
\end{minipage}
\\
\hline
\begin{tabular}{c}Convolution \\ Layer\end{tabular} &
\begin{minipage}{1.05\linewidth}
\begin{equation}
    {\boldsymbol\delta}_{\tau q}^{C}((s-1+l)h_{t};\boldsymbol{\omega}) \defeq \nabla_{{\boldsymbol{z}}_{\tau}^{C}((s-1+l)h_{t};\boldsymbol{\omega})} J_{\tau q}(\boldsymbol{\omega}) = h_{t}\sum\limits_{j=0}^{\min(s-1,q-l)}\alpha_{j}{\boldsymbol\delta}_{\tau q}^{+}((s-1+l+j;\boldsymbol{\omega})h_{t})\nonumber
\end{equation}
\end{minipage}
\\
\hline
\begin{tabular}{c}Reshape \\ Layer\end{tabular} &
\begin{minipage}{1.05\linewidth}
\begin{equation}
    \left[{\Delta}_{\tau q}^{R}((s-1+l)h_{t};\boldsymbol{\omega})\right]_{ij} \defeq \frac{\partial}{\partial \left[Z_{\tau}^{R}((s-1+l)h_{t};\boldsymbol{\omega})\right]_{ij}}J_{\tau q}(\boldsymbol{\omega}) =  \left[{\boldsymbol\delta}_{\tau q}^{C}((s-1+l)h_{t};\boldsymbol{\omega})\right]_{i}\left[-\frac{c}{h_{x}}\boldsymbol{\omega}_{1}+\frac{\nu}{h_{x}^{2}}\boldsymbol{\omega}_{2}\right]_{j}\nonumber
\end{equation}
\end{minipage}
\\
\hline
\begin{tabular}{c}Assignment \\ Layer\end{tabular} &
\begin{minipage}{1.05\linewidth}
\begin{equation}
    {\boldsymbol\delta}_{\tau q}^{A}((s-1+l)h_{t};\boldsymbol{\omega}) \defeq \nabla_{{\boldsymbol{z}}_{\tau}^{A}((s-1+l)h_{t};\boldsymbol{\omega})} J_{\tau q}(\boldsymbol{\omega})=\sum\limits_{j=1}^{n}R_{j}^{T}({\Delta}_{\tau q}^{R}((s-1+l)h_{t};\boldsymbol{\omega})\mathbf{e}_{j})\nonumber
\end{equation}
\end{minipage}
\\
\hline
\begin{tabular}{c}Input \\ Layer\end{tabular} &
\begin{minipage}{1.05\linewidth}
\begin{equation}
    {\boldsymbol\delta}_{\tau q}^{I}((s-1+l)h_{t};\boldsymbol{\omega}) \defeq \nabla_{{\boldsymbol{z}}_{\tau}^{I}((s-1+l)h_{t};\boldsymbol{\omega})} J_{\tau q}(\boldsymbol{\omega}) = {\boldsymbol\delta}_{\tau q}^{+}((s-1+l)h_{t};\boldsymbol{\omega})+U^{T}{\boldsymbol\delta}_{\tau q}^{A}((s-1+l)h_{t};\boldsymbol{\omega})\nonumber
\end{equation}
\end{minipage}
\\
\hline
\end{tabularx}
\end{center}
\caption{Intermediate quantities useful for backpropagation.  These values are computed for each $q=0,1,\ldots,Q-1$ and for $l=q,q-1,\ldots,0$.}
\label{tab:backprop1}
\end{table}

\begin{table}[p]
\begin{center}
\begin{tabularx}{1.05\linewidth}{|c|X|}
\hline
\begin{tabular}{c}Convolution \\ Layer\end{tabular} &
\begin{minipage}{1.05\linewidth}
\begin{equation}
    {\boldsymbol\delta}_{\tau q}^{C}(lh_{t};\boldsymbol{\omega}) \defeq \nabla_{{\boldsymbol{z}}_{\tau}^{C}(lh_{t};\boldsymbol{\omega})} J_{\tau q}(\boldsymbol{\omega}) = h_{t}\sum\limits_{j=s-l}^{\min(s,s+q-l)}\alpha_{j-1}{\boldsymbol\delta}_{\tau q}^{+}((l+j-1)h_{t};\boldsymbol{\omega})\nonumber
\end{equation}
\end{minipage}
\\
\hline
\end{tabularx}
\end{center}
\caption{Intermediate quantities useful for backpropagation. These values are computed for each $q=0,1,\ldots,Q-1$ and for $l=s-2,s-3,\ldots,0$.}
\label{tab:backprop2}
\end{table}
\end{landscape}

\section{Numerical Experiments} \label{sec:numexp}

Numerical experiments were performed to investigate the circumstances, given the many possible parameter choices detailed so far, for which the set of weights, $\boldsymbol{\omega}$, determined through the minimization of $J$ lead to stable solutions of \eqref{eq:ContinousProblem}. Since the decay rate of the coefficients of the Fourier series of the initial data have a correspondence to the differentiability of the initial data, three variations on initial conditions were also considered as a proxy for differing smoothness.  That is, for each of $p=0,2,8$, $T$ cases of initial data are generated via \eqref{eq:continuous_solution} with $b_{\eta}(0) = b_{-\eta}(0) = \frac{1}{\eta^{p}}\hat{b}_{\eta,T}$, for $\eta=0,1,\ldots,(N-1)/2$, and $b_{\eta}(0) = b_{-\eta}(0)=0$ otherwise.  Here, $\hat{b}_{\eta,T}$ is a value chosen from the random uniform distribution on $(-1,1)$, and solutions to \eqref{eq:ContinousProblem} are given by \eqref{eq:continuous_solution} for all values of $x\in\mathbb{R}$ and $t>0$.  

Aside from floating point errors, this training data is exact.  In a more complex setting where exact solutions to the PDE are not available, high-fidelity numerical solutions need to be used as training data.  These numerical solutions would have their own associated truncation errors and would be subject to any instabilities in the methods used to generate them.  The errors and instabilities would then have a significant influence on the set of weights generated by the procedure outlined herein if special care is not taken to mitigate their impact.

In this work $n$ is chosen to be an odd number, so that the local node set is centered on $x_{k}$, that is, $\mathcal{N}_{k,n} = \left\{x_{k}+(-(n-1)/2+j-1)h_{x}\right\}_{j=1}^{n}$.  Define $D'$ to be the matrix, with eigenvalues $\lambda'$, constructed from the second-order centered differences ($n=3$ with $\psi_{k,j}$ a Lagrange polynomial basis element in \eqref{eq:localfdapprox}) when approximating $\frac{\partial}{\partial x}$ and $\frac{\partial^{2}}{\partial x^{2}}$ with periodicity enforced.  That is, $D'$ is constructed as in \eqref{eq:Diffmatentries} from the weight vector $\boldsymbol{\omega}'=-c/h_{x}\mathbf{w}_{1}'+\nu/h_{x}^{2}\mathbf{w}_{2}'$, with
\begin{align}
\mathbf{w}_{1}'=\frac{1}{2}\left[\begin{array}{ccccc}\mathbf{0}_{(n-3)/2}^{T} & -1 & 0 & 1 & \mathbf{0}_{(n-3)/2}^{T}\end{array}\right]^{T}\nonumber
\end{align}
and
\begin{align}
\mathbf{w}_{2}'=\left[\begin{array}{ccccc}\mathbf{0}_{(n-3)/2}^{T} & 1 & 2 & 1 & \mathbf{0}_{(n-3)/2}^{T}\end{array}\right]^{T}.\nonumber
\end{align}

\begin{figure}
    \begin{center}
    \includegraphics[width=\textwidth]{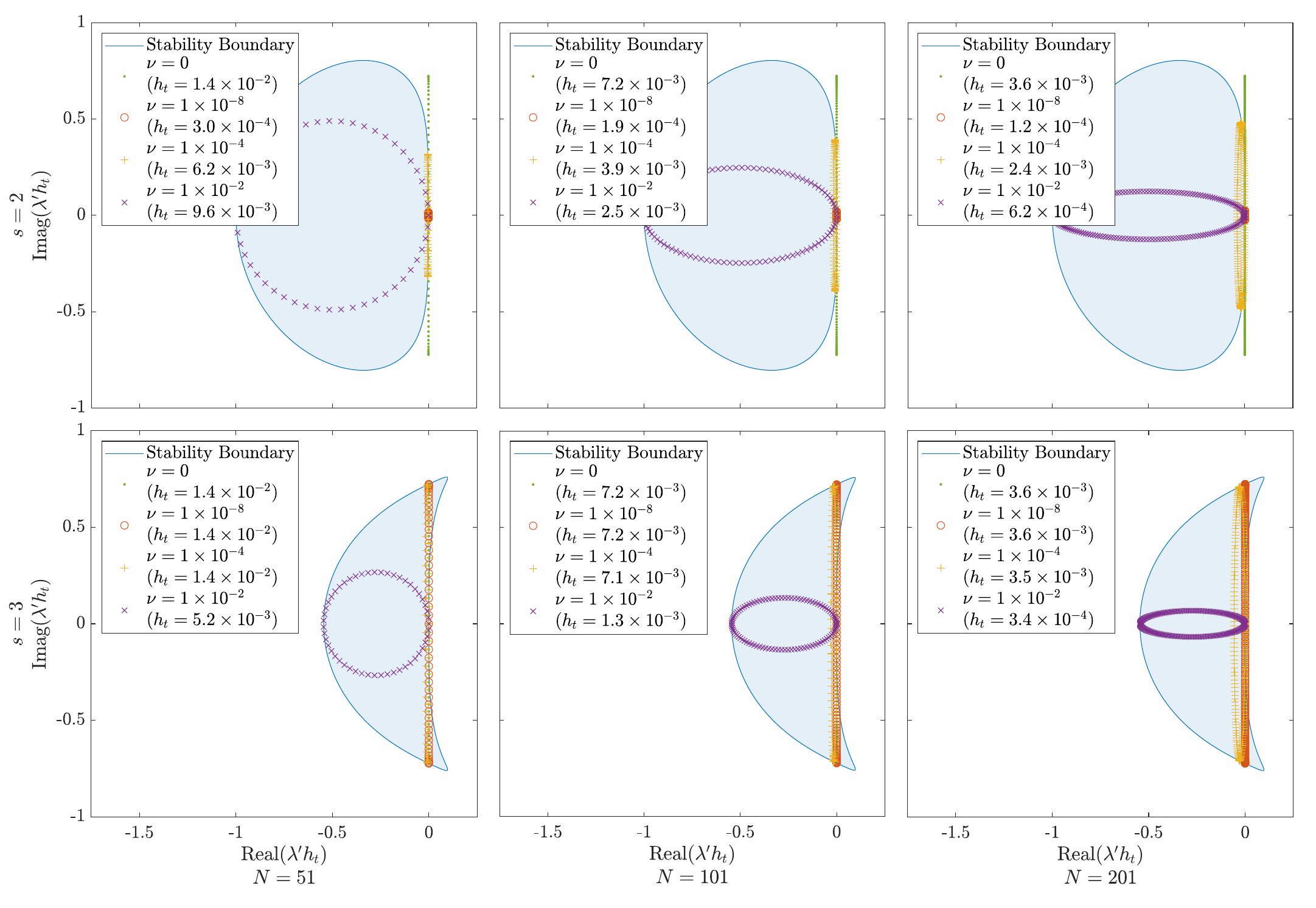}
    \end{center}
    \caption{Eigenvalues of $D'$ for various values of $\nu$ with $h_{t}=h_{t,N,\nu,s}$.  The value of $h_{t,N,\nu,s}$ is indicated in parentheses for each choice of $N$, $\nu$ and $s$.  The stability region (shaded) of AB-$s$ is shown with its stability boundary outlined by a solid curve.  Markers indicate the locations of the scaled eigenvalues. Notice that these choices of $h_{t}$ are (approximately) the largest such that the scaled eigenvalues fall inside the stability region.}
    \label{fig:StabRegionEvals}
\end{figure}

When $\nu=0$ and $c>0$ this matrix is antisymmetric, and so has purely imaginary eigenvalues.  Alternatively, when $\nu>0$ and $c=0$ the matrix $D'$ is symmetric with purely real eigenvalues.  Otherwise, the eigenvalues exist on a curve in the left half of the complex plane, with negative real part.  Further, the stability region of AB-$2$ is a closed curve contained entirely in the left half of the complex plane (along with the origin).  Therefore, when $\nu=0$ there is no nonzero value of $h_{t}$ such that $\lambda'h_{t}$ is contained in the stability region of the multi-step method.  For all other choices of $\nu>0$, there is a maximum value, denoted $h_{t,N,\nu,s}$, such that $\lambda'h_{t,N,\nu,s}$ is contained in the stability region of AB-$s$.  The value of $h_{t,N,\nu,s}$ is chosen numerically and iteratively utilizing a bisection process that brackets this critical value of $h_{t}$ and tests which half of the bracket produces eigenvalues contained entirely in the stability region.  In the case of $\nu=0$ and $s=2$, since there is no such value, $h_{t,N,0,2}$ is set to the value of $h_{t,N,0,3}$. Figure \ref{fig:StabRegionEvals} illustrates the stability regions and $\lambda'h_{t}$ when $h_{t}$ is set to these critical values for all choices of $N$, $s$ and $\nu$ investigated here.  The solid curve, known as the stability boundary, illustrates the boundary of the (shaded) stability region, while the markers (dependent on $\nu$) illustrate the locations of the scaled eigenvalues.

Experiments are performed for values of $h_{t}$ near these critical values of $h_{t,N,\nu,s}$, since a goal here is to determine whether minimizing $J$ can deliver a stable method for values of $N$, $h_{t}$, $\nu$ and $s$ where traditional methods may fail.   Table \ref{tab:Parameters} includes all parameter values used in this study, and each possible combination (for a total of $36450$) computational experiments were performed.  Further, in all cases the vector of weights is initialized as
\begin{align}
\boldsymbol{\omega}^{(0)} = \left[\begin{array}{cc}(\mathbf{w}_{1}')^{T} & (\mathbf{w}_{2}')^{T}\end{array}\right]^{T}\nonumber
\end{align}
for the quasi-Newton method.  Other choices of initial weight vectors are certainly possible, but this choice provided the most success in constructing stable methods for the combinations of parameters listed in table \ref{tab:Parameters}.

\begin{table}
\begin{tabularx}{\linewidth}{|X|X|}
\hline
Parameter Name (Description) & Set of Values \\
\hline
$c$ (advection coefficient/speed) & $1$ \\
\hline
$\nu$ (diffusion coefficient) &  $0$, $1\times10^{-4}$, and $1\times10^{-2}$ \\
\hline
$P$ (period) & $1$\\
\hline
$p$ (decay rate of Fourier series coefficients of initial data) & $2,4,8$ \\
\hline
$N$ (number of computational nodes in $S$ & $51,101,201$ \\
\hline
$n$ (number of nearest neighbors) & $3,5,7,9,11$ \\
\hline
$s$ (steps in multi-step method) & $2$ and $3$ \\
\hline
$h_{t}$ (time step in numerical integration) & $h_{t,N,\nu,s}$, $1.01h_{t,N,\nu,s}$, $1.1h_{t,N,\nu,s}$ \\
\hline
$Q$ (number of recurrent steps in loss) & $1,3,4,5,9$ \\
\hline
$T$ (number of training cases & $1,10,100$ \\
\hline
$\kappa_{\mbox{max}}$ (maximum number of iterations in quasi-Newton method) & $10,100,1000$ \\
\hline
\end{tabularx}
\caption{Parameter sets used in computational experiments.  These are the parameters that can be adjusted in the implementation of the method described in this paper.}
\label{tab:Parameters}
\end{table}

\begin{figure}
\begin{center}
\includegraphics[width=\textwidth]{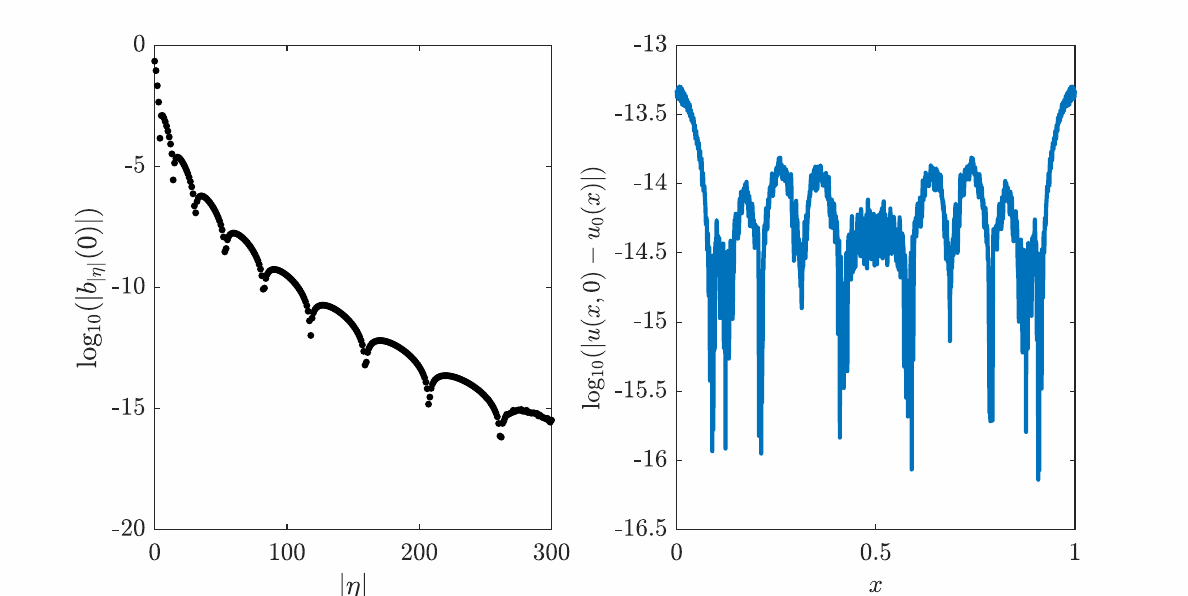}
\end{center}
\caption{Left: modulus of $b_{\eta}(0)$ for the bump function \eqref{eq:bumpfunction}. Right: log base 10 of the error in the approximation \eqref{eq:partialsum} of \eqref{eq:bumpfunction}.}\label{fig:BumpFunctionErrorandCoefficients}
\end{figure}

Once weights are computed by minimizing $J$, they are utilized to approximate the solution to \eqref{eq:ContinousProblem} for $t\in(0,20]$ and $x\in[0,1)$ beginning with initial data defined from the partial sum
\begin{align}
    u(x,t) = \sum\limits_{\eta=-300}^{300}b_{\eta}(0)\mbox{e}^{i\frac{2\pi \eta}{P}(x-ct)-\left(\frac{2\pi \eta}{P}\right)^{2}\nu t},\label{eq:partialsum}
\end{align}
with $t=0$ and $b_{\eta}(0)$ chosen to be the coefficients of the 1-periodic ($P=1$) extension of the bump function,
\begin{align}
    u_{0}(x) = \left\{\begin{array}{cc}\mbox{e}^{\frac{1}{(2(x-\lfloor x\rfloor)-1)^{2}-1}}, & x\notin\mathbb{Z}\\0, & x\in\mathbb{Z}\end{array}\right..\label{eq:bumpfunction}
\end{align}
These values of $b_{\eta}(0)$ are computed numerically from \eqref{eq:seriescoefficients} using Matlab's (R2022b) integral function with \verb+AbsTol+ set to \verb+eps+.  The modulus of $b_{\eta}(0)$ decreases with the modulus of $\eta$, although not monotonically, and reaches machine precision (near $10^{-16}$) when $\lvert\eta\rvert\approx300$.   This coincides with the initial data being continuous and smooth in the sense that all of its derivatives are continuous as well.  An illustration of these coefficients and the absolute difference between $u(x,0)$ from \eqref{eq:partialsum} and $u_{0}(x)$ from \eqref{eq:bumpfunction} is shown in figure \ref{fig:BumpFunctionErrorandCoefficients}.  Notice that for this initial data, exact solutions to \eqref{eq:ContinousProblem} are given by \eqref{eq:partialsum}.  Therefore, to evaluate the forward error in the approximate solution, the vector $\mathbf{u}(t)$ is constructed at discrete values of $t$, equally spaced by $h_{t}$, defined entrywise as $[\mathbf{u}(lh_{t})]_{k}=u(x_{k},t)$ with $u$ as in \eqref{eq:partialsum} and for $l\geq0$.  Further, to kickstart the AB-$s$ method \eqref{eq:ms_method}, it is useful to set $\tilde{\mathbf{u}}(lh_{t})=\mathbf{u}(lh_{t})$, $l=0,1,\ldots,s-1$.  It should also be noted that this initial data is not contained in the training data for any of $p=0,2,8$.

\section{Key Observation on Recurrent Loss and Stability} \label{sec:keyobs}

The experiments outlined in section \ref{sec:numexp} are intended, in part, to assist in understanding the mechanism by which recurrent loss can promote stability in the numerical solution of a PDE.  Remarkably, without imposing any explicit conditions on the spectrum of the matrix $D$, minimizing $J$ results in a set of eigenvalues such that $\{\lambda_{\eta}h_{t}\}_{\eta=-\lfloor N/2\rfloor}^{\lfloor N/2\rfloor}$ (with $\lambda_{\eta}$ as in \eqref{eq:Deigenvalues}) adapts, to some extent, to the stability region of the multi-step method under consideration.  In some cases, the weights in $D$ are altered in such a way that the set of scaled eigenvalues falls entirely within the stability region, resulting in a stable method.  The dependence of the adaptation of the scaled eigenvalues to the stability region (and stability characteristics of the resulting numerical method) on the various parameters discussed in section \ref{sec:numexp} is complicated, with no discernible patterns arising.  The remainder of this section discusses observations relative to $Q$, $T$, $\kappa_{\mbox{max}}$, $p$, and $N$. Note that a set of weights for approximating the action of $\mathcal{L}$ for which the overall numerical method is stable when paired with AB-$s$ will be denoted a \textit{stable set of weights} for the remainder of this study.

\begin{figure}
\begin{center}
\includegraphics[width=\textwidth]{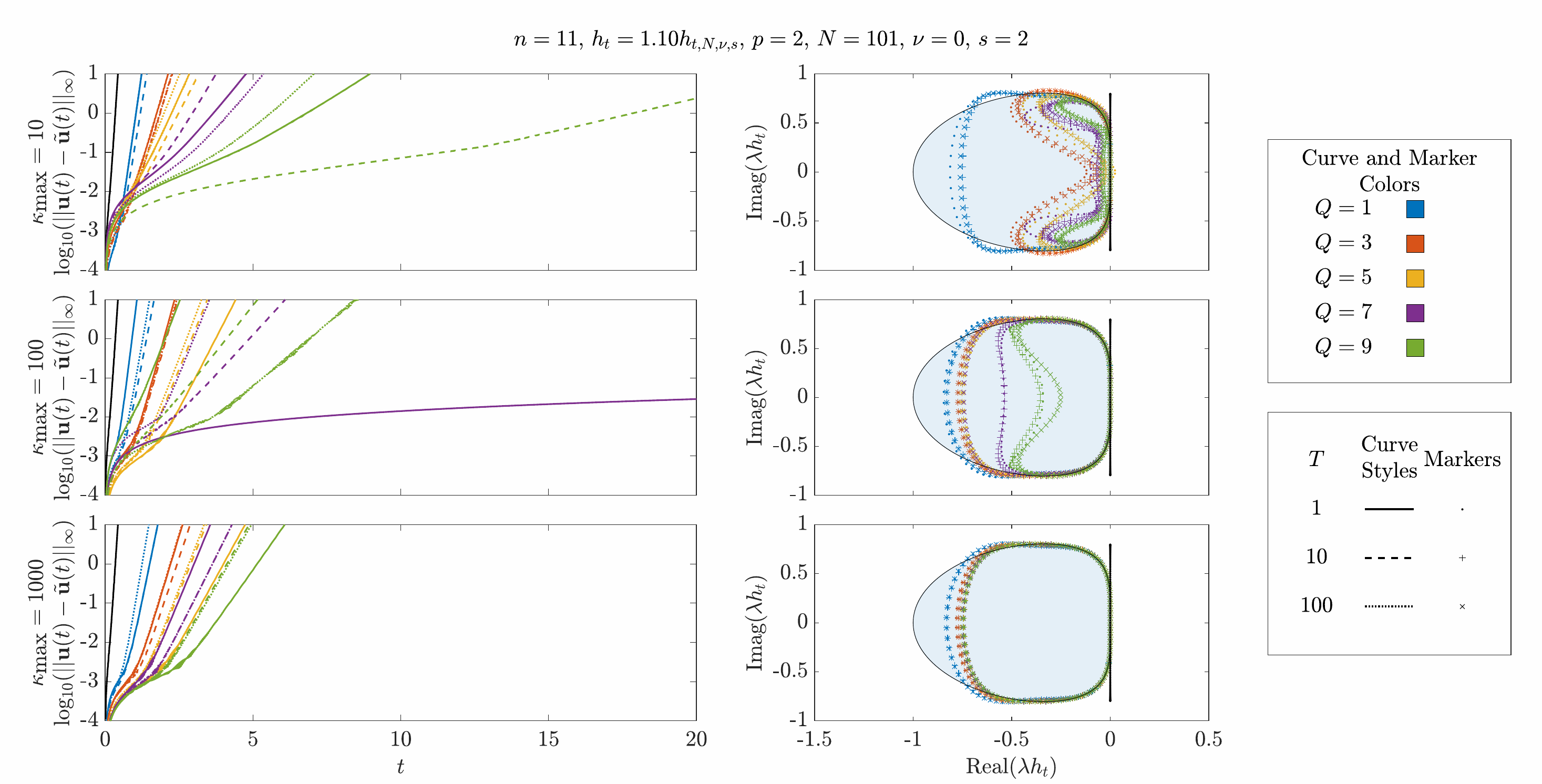}
\end{center}
\caption{Left column: log base 10 of the infinity norm absolute forward error in the numerical solution computed using the weights generated by minimizing $J$.  The colors indicate varying values of $Q$ and the curve styles depict different values of $T$.  The remaining parameter values are given in the figure title. Right column: eigenvalues of $D$, constructed from these weights, with markers this time indicating different values of $T$.  Rows correspond to different values of $\kappa_{\mbox{max}}$, increasing from top to bottom.  In all of the rows, the black curves in the left column and the black dots in the right column correspond to the absolute forward error and eigenvalues, respectively, of the second order centered difference method (with corresponding matrix $D'$). Parameter choices are indicated in the title of the plot.  Notice that a single set of weights is determined for $\kappa_{\mbox{max}}=100$, $Q=7$ and $T=1$, but for not other set of these three parameters depicted here.} \label{fig:ErrorandEigenvalues_s_2_n_11_neglog10nu_Inf_N_101_htmult_1_10_p_2}
\end{figure}

Consider figure \ref{fig:ErrorandEigenvalues_s_2_n_11_neglog10nu_Inf_N_101_htmult_1_10_p_2}, with the rows of plots corresponding to different values of $\kappa_{\mbox{max}}$, increasing from top to bottom.  The left column of subplots displays log base 10 of the infinity norm absolute forward error in the numerical solution computed using the weights generated by minimizing $J$.  The colors indicate varying values of $Q$ and the curve styles depict different values of $T$.  The remaining parameter values are given in the figure title.  Similarly, the right column of subplots depicts the eigenvalues of $D$, constructed from these weights, with markers this time indicating different values of $T$.  In all of the rows, the black curves in the left column and the black dots in the right column correspond to the absolute forward error and scaled eigenvalues, respectively, of the second order centered difference method (with corresponding matrix $D'$).  Here $\nu=0$ and $s=2$ so that second order centered difference method is unstable when paired with AB-2 for any value of $h_{t}>0$.  This instability produces solutions with errors that increase exponentially, even for $t\ll1$, with a rate dependent on the size of certain eigenvalues for which $\lambda_{\eta}h_{t}$ falls outside of the stability region.  However, even for $h_{t}=1.10h_{t,N,\nu,s}$ a stable set of weights can be determined through the minimization of $J$.  This can be seen in the second row of subplots when $Q=7$ and $T=1$.  Over $\kappa_{\mbox{max}}=100$ iterations of the quasi-Newton method, a set of weights is determined for which the forward error in the solution remains well below $10^{-1}$ for all $lh_{t}\in[0,20]$, $l>0$.  Referring to the corresponding sets of eigenvalues in the right subplots, the set of parameters corresponding to this stable solution also produces a set of eigenvalues so that $\lambda_{\eta}h_{t}$ falls entirely inside the stability region.

\begin{figure}
\begin{center}
\includegraphics[width=\textwidth]{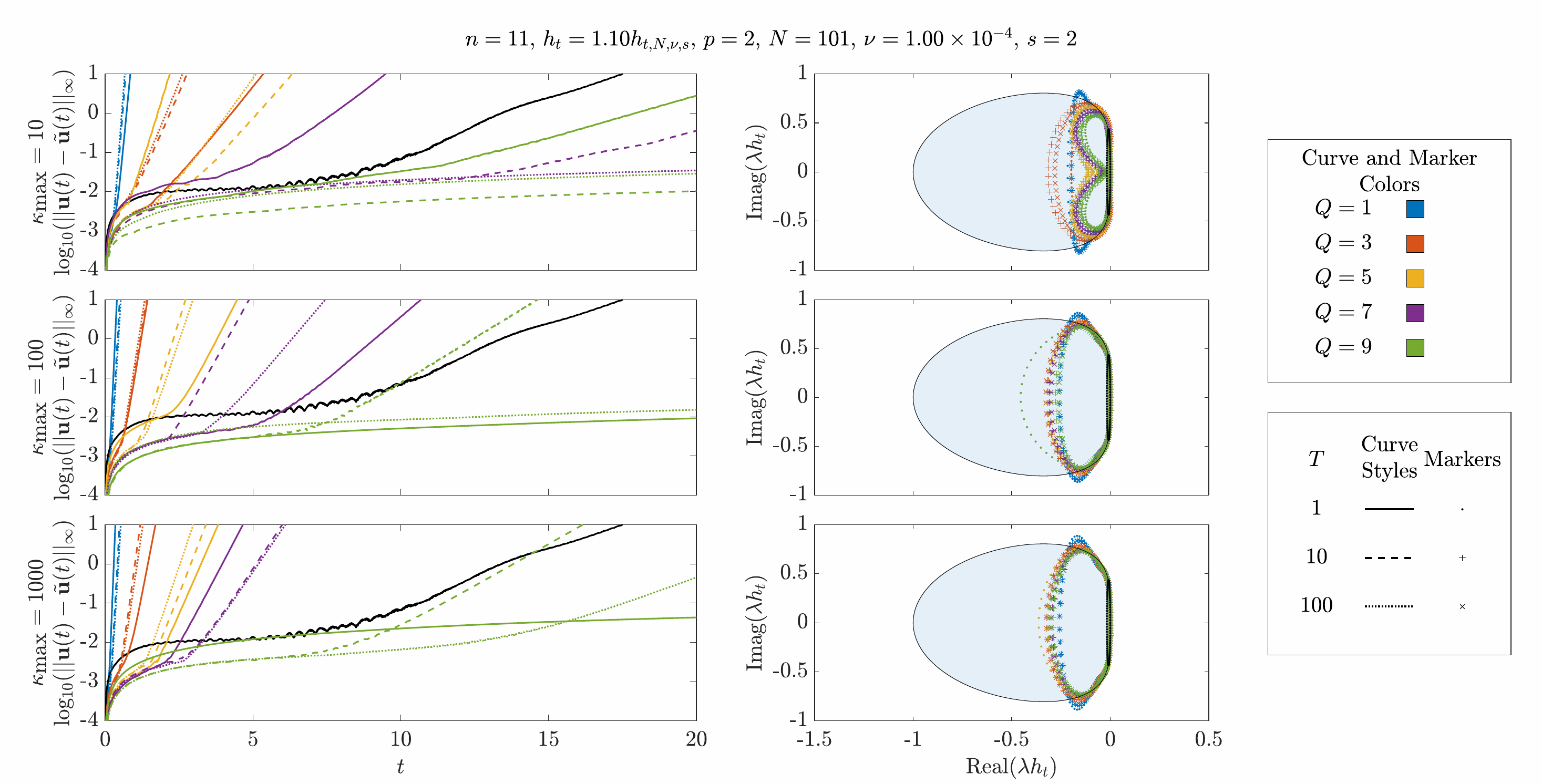}
\end{center}
\caption{See figure \ref{fig:ErrorandEigenvalues_s_2_n_11_neglog10nu_Inf_N_101_htmult_1_10_p_2} for a discussion of what is depicted in each subplot. Parameter choices are indicated in the title of the plot. Notice that more stable sets of weights are determined than the case depicted in \ref{fig:ErrorandEigenvalues_s_2_n_11_neglog10nu_Inf_N_101_htmult_1_10_p_2}.}\label{fig:ErrorandEigenvalues_s_2_n_11_neglog10nu_004_N_101_htmult_1_10_p_2}
\end{figure}

Changing $\nu$ from $0$ to $10^{-4}$ leads to a case where, for each $N$, there is a maximum value $h_{t}=h_{t,N,\nu,s}$ for which the numerical method pairing the second order centered difference method and AB-2 is stable.  The right subplots of figure \ref{fig:ErrorandEigenvalues_s_2_n_11_neglog10nu_004_N_101_htmult_1_10_p_2} illustrate that, for the depicted set of parameters, the values of $\lambda_{\eta}'h_{t}$ are located on a curve that falls partially within the stability region of AB-2.  Further, figure \ref{fig:StabRegionEvals} illustrates that when $s=2$ the maximum distance of any value of $\lambda_{\eta}'h_{t}$ from the stability domain when $\nu=10^{-4}$ is less than when $\nu=0$.  So, leaving all other parameters the same as those used in the results presented in figure \ref{fig:ErrorandEigenvalues_s_2_n_11_neglog10nu_Inf_N_101_htmult_1_10_p_2}, figure \ref{fig:ErrorandEigenvalues_s_2_n_11_neglog10nu_004_N_101_htmult_1_10_p_2} illustrates that the minimization of $J$ now allows the determination of multiple stable sets of weights.

\begin{figure}
\begin{center}
\includegraphics[width=\textwidth]{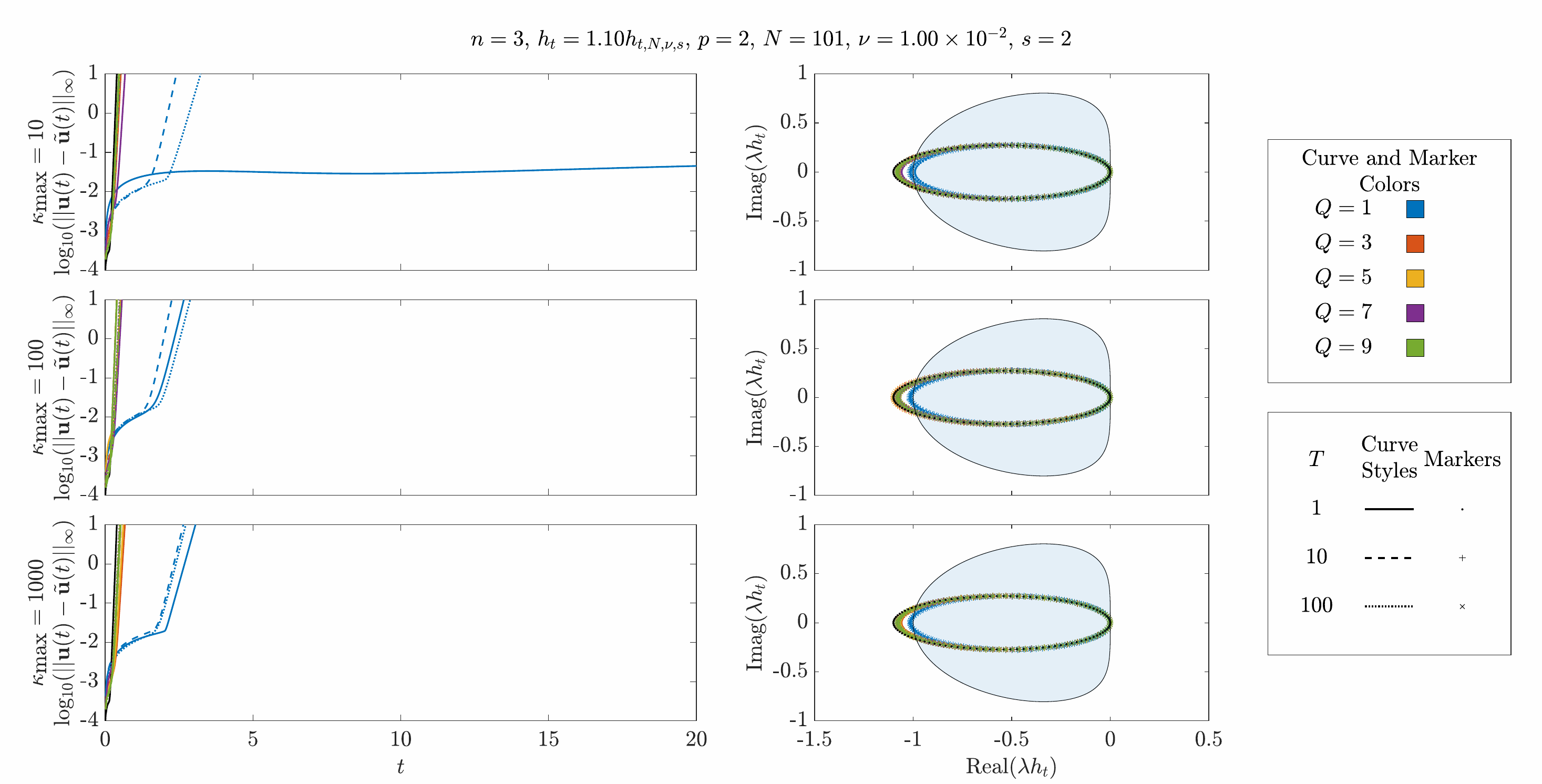}
\end{center}
\caption{See figure \ref{fig:ErrorandEigenvalues_s_2_n_11_neglog10nu_Inf_N_101_htmult_1_10_p_2} for a discussion of what is depicted in each subplot. Parameter choices are indicated in the title of the plot. Notice that a stable set of weights is determined only for $\kappa_{\mbox{max}}=10$, $Q=1$ and $T=1$.  That is, in this case inclusion of multiple time steps in the recurrent loss does not improve stability.}\label{fig:ErrorandEigenvalues_s_2_n_3_neglog10nu_002_N_101_htmult_1_10_p_2}
\end{figure}

Noticeably, in figures \ref{fig:ErrorandEigenvalues_s_2_n_11_neglog10nu_Inf_N_101_htmult_1_10_p_2} and \ref{fig:ErrorandEigenvalues_s_2_n_11_neglog10nu_004_N_101_htmult_1_10_p_2} the values of $t$ for which the growth in $\lVert \mathbf{u}(t)-\tilde{\mathbf{u}}(t)\rVert_{\infty}$ transitions from polynomial to exponential increases with $Q$ in most cases.  That is, the number of time steps included in the recurrent loss function can have a positive impact on generating weights for approximating the action of $\mathcal{L}$ that lead to a stable numerical method.  Still, the choice of $Q$ that will produce a stable set of weights is often unpredictable.  For instance, in figure \ref{fig:ErrorandEigenvalues_s_2_n_11_neglog10nu_Inf_N_101_htmult_1_10_p_2}, the choice of $Q=7$ produces the only stable set of weights, even when $Q=9$ was also considered.  Figure \ref{fig:ErrorandEigenvalues_s_2_n_3_neglog10nu_002_N_101_htmult_1_10_p_2} further illustrates that for certain parameter sets, increases in $Q$ can have little to no impact on the ability to determine a stable set of weights.  With $n$ changed from $11$ to $3$ and $\nu$ to $10^{-2}$, the only stable set of weights is found when $\kappa_{\mbox{max}}=10$, $Q=1$, and $T=1$.  Figures \ref{fig:ErrorandEigenvalues_s_2_n_11_neglog10nu_Inf_N_101_htmult_1_10_p_2}, \ref{fig:ErrorandEigenvalues_s_2_n_11_neglog10nu_004_N_101_htmult_1_10_p_2} and \ref{fig:ErrorandEigenvalues_s_2_n_3_neglog10nu_002_N_101_htmult_1_10_p_2} illustrate some of the most extreme observed behavior relative to the choice of $Q$; however, they demonstrate an important point that is true for nearly all computed parameter sets.  That is, increases in $Q$ do not reliably correspond to an increased ability to generate a stable set of weights. This unpredictability with respect to $Q$ is likely due to a combination of factors.  First, increasing $Q$ leads to a larger number of layers that must be traversed in the process of computing the gradients via backpropagation, escalating the chances that the well known problem of vanishing gradients will limit the effectiveness of the training procedure.  Second, there is a competition between the goals of accuracy and stability, but with stability only encouraged by attempting to minimize forward error over more than a single step.  A local minimizer of the loss function is then not guaranteed to be a stable set of weights but one that achieves the smallest value of the loss out of nearby weight sets.  

Likewise, as $T$ increases the set of training data includes a larger subset of the set of $1$-periodic functions.  Therefore, it might be expected that the computed weights converge to a particular set of values.  However, there are clearly combinations of $T$ and $Q$ for which the $T=1$ and $T=100$ cases perform more alike than the $T=10$ case (see, for instance the curves for $Q=9$ in the second row of subplots in figure \ref{fig:ErrorandEigenvalues_s_2_n_11_neglog10nu_004_N_101_htmult_1_10_p_2}).  Now, it is likely that the choices of $T$ explored here are not large enough to demonstrate such convergence; however, if the convergence is to a set of weights for which the method for solving \eqref{eq:ContinousProblem} is unstable or low accuracy, then there is no benefit to consider even larger values of $T$.

It should be expected that the weight set computed using the quasi-Newton method converges superlinearly to a stationary point of $J$ (see, e.g., \cite{Nocedal_Wright_2006}).  The rows of subplots in each of figures \ref{fig:ErrorandEigenvalues_s_2_n_11_neglog10nu_Inf_N_101_htmult_1_10_p_2}, \ref{fig:ErrorandEigenvalues_s_2_n_11_neglog10nu_004_N_101_htmult_1_10_p_2} and \ref{fig:ErrorandEigenvalues_s_2_n_3_neglog10nu_002_N_101_htmult_1_10_p_2} demonstrate that, while minimization of the forward error in $\tilde{\mathbf{u}}(lh_{t})$ for $l=1,2,\ldots,Q$ simultaneously can promote both stability and improved accuracy, the number of iterations $\kappa_{\mbox{max}}$ required to achieve these improvements is uncertain.  Ideally, once a set stable set of weights is found for a particular choice of $Q$ and $T$, additional iterations in the minimization of $J$ should attempt to reduce the forward error, particularly over the set of $Q$ steps under consideration in the recurrent loss.  This appears to be the case in figure \ref{fig:ErrorandEigenvalues_s_2_n_11_neglog10nu_004_N_101_htmult_1_10_p_2} when $\kappa_{\mbox{max}}=10,100$ for $Q=9$ and $T=100$; however, it is not so for $\kappa_{\mbox{max}}=100,1000$ when $Q=9$ and $T=1$, although this may be attributable to the small number of training cases. Further, in figure \ref{fig:ErrorandEigenvalues_s_2_n_3_neglog10nu_002_N_101_htmult_1_10_p_2} a stable set of weights is found only when $\kappa_{\mbox{max}}=10$ and increases to 100 and 1000 no longer produce stable sets of weights.   Similarly, in figure \ref{fig:ErrorandEigenvalues_s_2_n_11_neglog10nu_Inf_N_101_htmult_1_10_p_2} the stable set of weights is determined only when $\kappa_{\mbox{max}}=100$.

\begin{figure}
\begin{center}
\includegraphics[width=\textwidth]{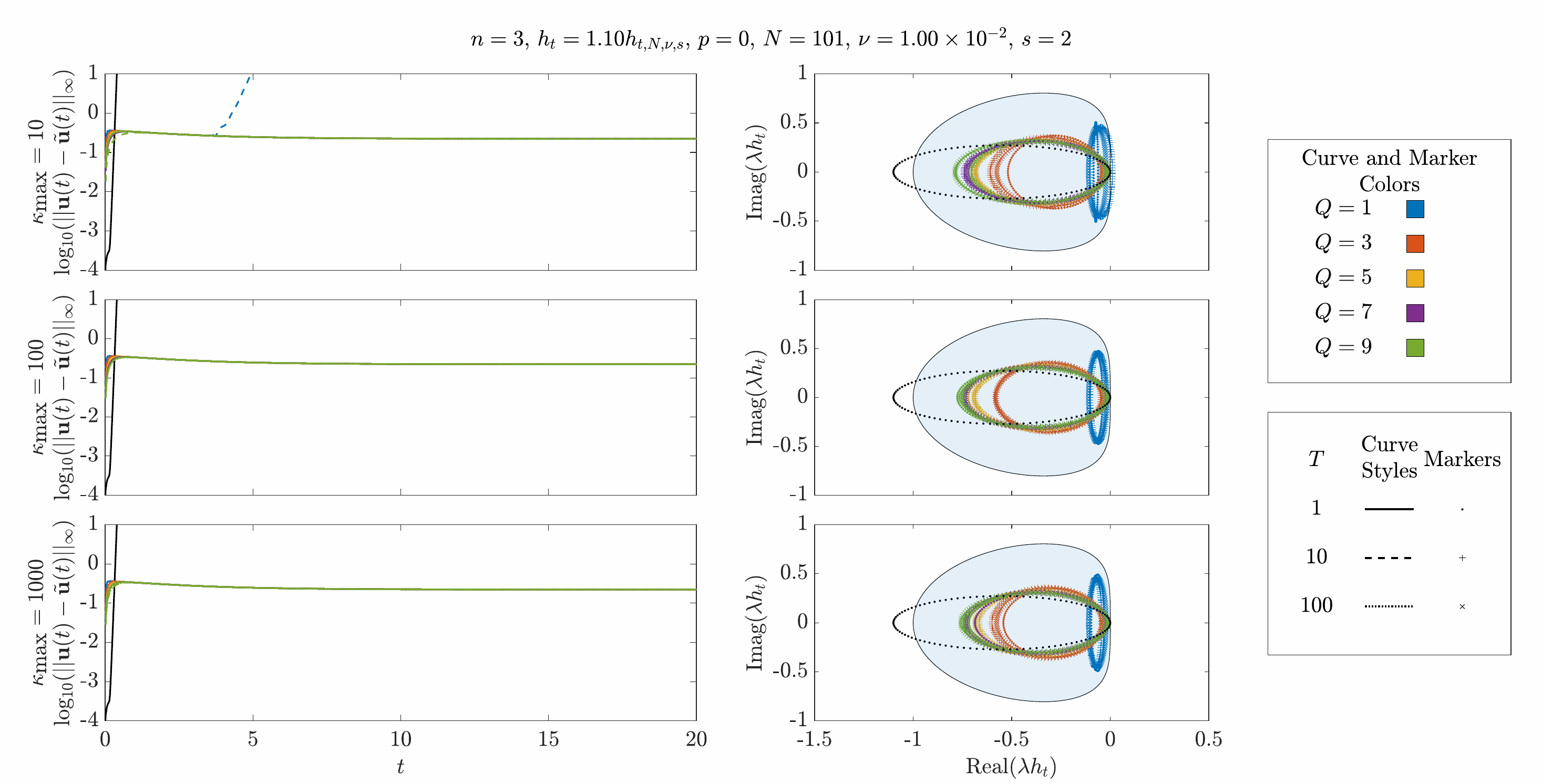}
\end{center}
\caption{See figure \ref{fig:ErrorandEigenvalues_s_2_n_11_neglog10nu_Inf_N_101_htmult_1_10_p_2} for a discussion of what is depicted in each subplot. Parameter choices are indicated in the title of the plot. Notice that training on data with $p=0$ produces more stable sets of weights since more of the values of $\lambda_{\eta}h_{t}$ are emphasized.}\label{fig:ErrorandEigenvalues_s_2_n_3_neglog10nu_002_N_101_htmult_1_10_p_0}
\end{figure}

The value of $p$, which impacts the behavior of the initial data (that is, $b_{|\eta|}(0) = O(\eta^{-p}$) as $\eta\to\infty$), is chosen to mimic varying degrees of smoothness of the initial data. As $p$ increases the contribution of $\mathbf{v}_{|\eta|}$, $|\eta|>1$, in \eqref{eq:EigenvectorSolution} diminishes.  Further, the eigenvalues \eqref{eq:Deigenvalues} corresponding to larger $|\eta|$ typically have modulus that is larger than those for small $|\eta|$.  That is, $\lambda_{\eta}h_{t}$ is more likely to fall outside of the stability region of the AB-$s$ method for larger $|\eta|$ than for small.  Therefore, the diminished influence of these eigenvectors discounts the influence of the corresponding eigenvalues on instability, particularly for the small values of $t$ considered in the process of minimizing $J$.  Practically speaking, this has translated to fewer combinations of $s$, $n$, $\nu$, $N$ and $h_{t}$ for which any weight sets leading to a stable method being generated as $p$ increases since the objective detects less growth from the scaled eigenvalues outside the stability region.  This is clear when comparing figure \ref{fig:ErrorandEigenvalues_s_2_n_3_neglog10nu_002_N_101_htmult_1_10_p_0} to figure \ref{fig:ErrorandEigenvalues_s_2_n_3_neglog10nu_002_N_101_htmult_1_10_p_2} where the only parameter that has changed is the value of $p$.  With an increased influence of eigenvalues corresponding to larger values of $\lvert\eta\rvert$ when $p=0$, more stable sets of weights are determined as shown in figure \ref{fig:ErrorandEigenvalues_s_2_n_3_neglog10nu_002_N_101_htmult_1_10_p_0} (surprisingly for all but one pair of $Q$ and $T$).  Still, figure \ref{fig:ErrorandEigenvalues_s_2_n_3_neglog10nu_002_N_101_htmult_1_10_p_0} also demonstrates what can be lost when choosing $p$ too low.  Comparing again to figure \ref{fig:ErrorandEigenvalues_s_2_n_3_neglog10nu_002_N_101_htmult_1_10_p_2}, when the eigenvalues that fall outside the stability region are given more emphasis in the case of $p=0$ the accuracies of the resulting numerical methods appear to suffer in comparison to the $p=2$ case since the minimization of $J$ must contend with the instability to a greater extent, impacting further improvements in accuracy.

\begin{figure}
\begin{center}
\includegraphics[width=\textwidth]{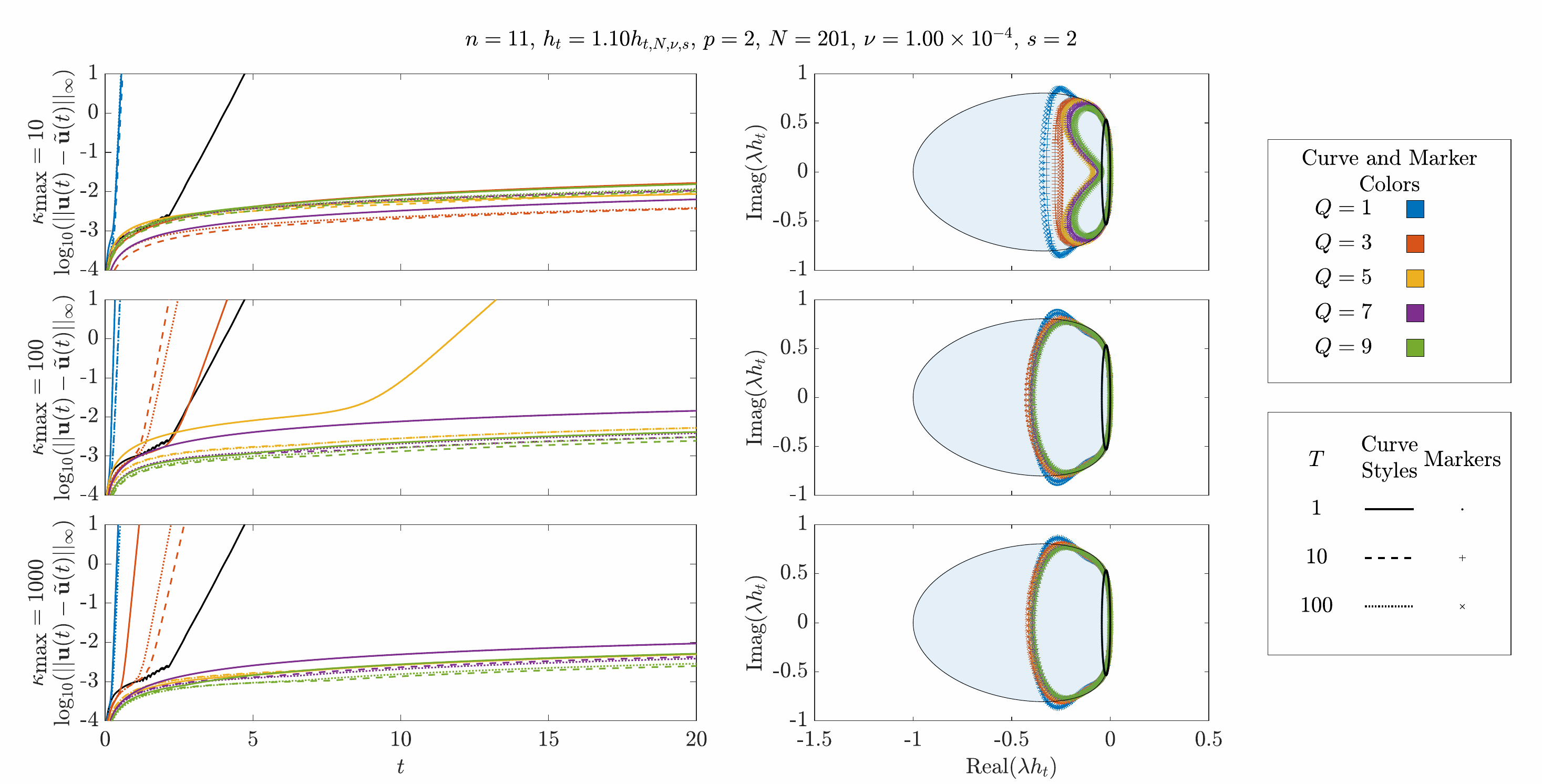}
\end{center}
\caption{See figure \ref{fig:ErrorandEigenvalues_s_2_n_11_neglog10nu_Inf_N_101_htmult_1_10_p_2} for a discussion of what is depicted in each subplot. Parameter choices are indicated in the title of the plot. Notice that an increase in $N$ can improve the accuracy in the solution.}\label{fig:ErrorandEigenvalues_s_2_n_11_neglog10nu_004_N_201_htmult_1_10_p_2}
\end{figure}

Figure \ref{fig:ErrorandEigenvalues_s_2_n_11_neglog10nu_004_N_201_htmult_1_10_p_2} is generated from the same parameter set as figure \ref{fig:ErrorandEigenvalues_s_2_n_11_neglog10nu_004_N_101_htmult_1_10_p_2} except that $N$ is increased from 101 to 201. Notice that in this case, the choices of $Q$ and $T$ for which a stable set of weights are found when $N=101$ also produce a stable set of weights when $N=201$.  These stable sets of weights appear to be more accurate for the larger value of $N$, particularly for $\kappa_{\mbox{max}}=100,1000$.  Increases in the node set size, $N$, lead to smaller values of $h_{t,N,\nu,s}$ required for the centered difference method to be stable (excluding the case of $s=2$ and $\nu=0$); however, the centered difference is also more accurate for larger $N$. Figures \ref{fig:ErrorandEigenvalues_s_2_n_11_neglog10nu_004_N_101_htmult_1_10_p_2} and \ref{fig:ErrorandEigenvalues_s_2_n_11_neglog10nu_004_N_201_htmult_1_10_p_2} hint that the same behavior may be possible for numerical methods utilizing the weight sets constructed by minimizing $J$. Further, there are many pairs of $Q$ and $T$ for which stable sets of weights are found when $N=201$ but not for $N=101$ and the same value of $\kappa_{\mbox{max}}$.  Unfortunately, except in rare cases like the one depicted here, increases in $N$ do not reliably lead to the determination of weights that promote stability or greater accuracy in the resulting numerical methods.

\subsection{Truncation Error}

One important concern with any method for approximating the solution to a PDE is how the error in the solution scales under refinement of the spatial discretization or of the size of the steps in the time (or possibly time-like) variable. Here the truncation error of the time-stepping scheme, AB-$s$, is fixed, so some computational experiments were performed to estimate the truncation errors in the approximate action of the operators acting with respect to the spatial variables..  To estimate the local truncation error, the value of $\lVert \mathbf{u}(sh_{t})-\tilde{\mathbf{u}}(sh_{t})\rVert_{\infty}$ is computed, given \eqref{eq:partialsum} and \eqref{eq:bumpfunction}, for the numerical methods that utilize the weights generated for $n=11$,  $h_{t}=1.10h_{t,N,\nu,s}$, $p=2$, $N=201$, $\nu=1.00\times10^{-4}$ and $s=2$ but now with $h_{t}=10^{-5}$ and $N\in[11,10001]$.  That is, first, the weights for approximating the action of $\mathcal{L}$ are generated using the same set of parameters as in figure \ref{fig:ErrorandEigenvalues_s_2_n_11_neglog10nu_004_N_201_htmult_1_10_p_2}, then the numerical method that pairs AB-2 with these weights is used to compute solutions after a single step in $t$ for various values of $N$ (refining the spatial discretization) and fixing $h_{t}$.  The choice of $h_{t}$ is less than the smallest computed value of $h_{t,N,\nu,s}$ for these choices of $N$ to avoid any unnecessary growth in the solution, even after a single step, for stable methods.

\begin{figure}
\begin{center}
\includegraphics[width=\textwidth]{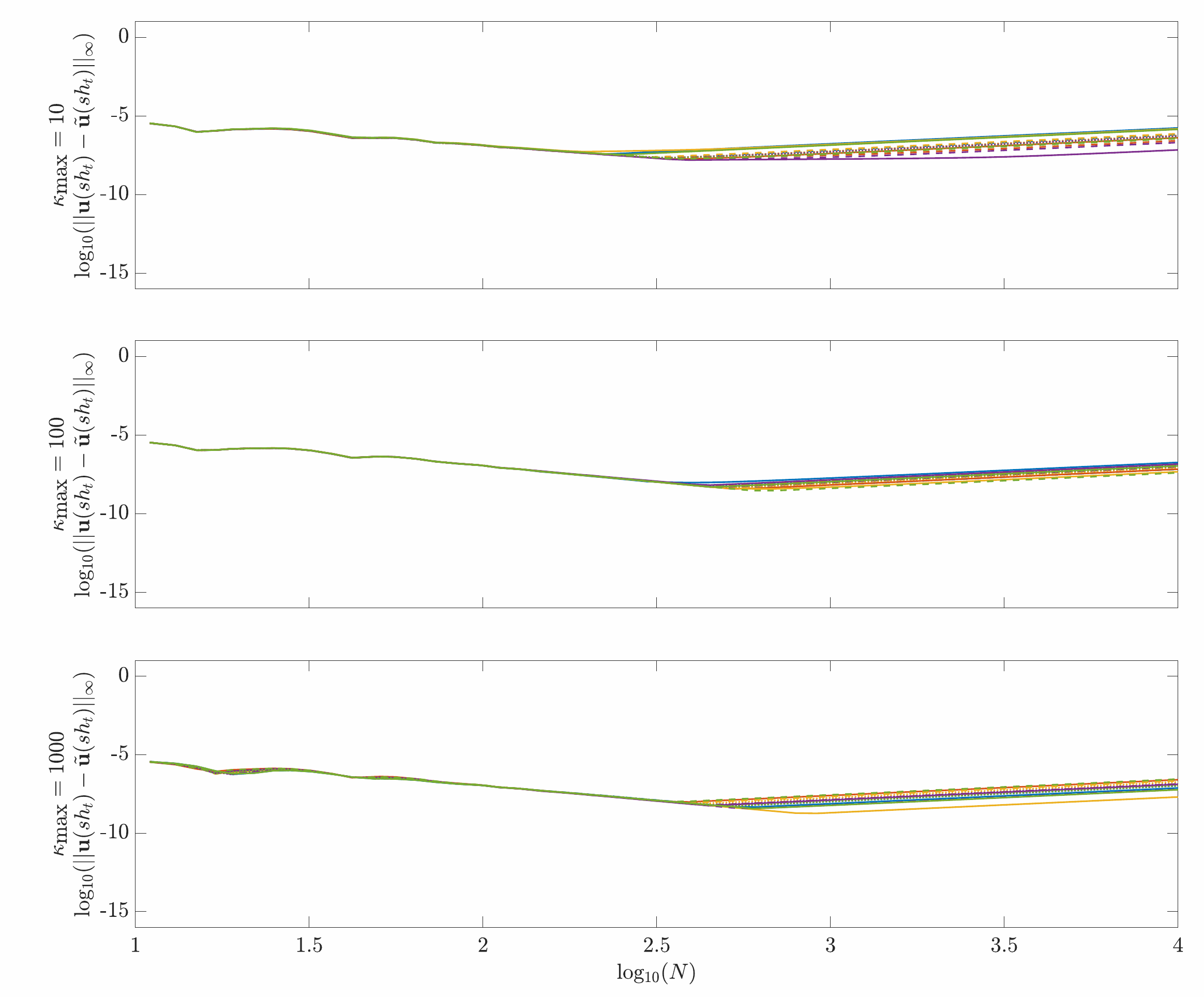}
\end{center}
\caption{Error when taking a single step time step using AB-2 and the weights generated with parameters $n=11$,  $h_{t}=1.10h_{t,N,\nu,s}$, $p=2$, $N=201$, $\nu=1.00\times10^{-4}$ and $s=2$. These parameters correspond to those in figure \ref{fig:ErrorandEigenvalues_s_2_n_11_neglog10nu_004_N_201_htmult_1_10_p_2}.  Here $h_{t}$ is fixed and the value of $N$ is allowed to range over the interval $[11,10001]$, indicating refinement of the spatial discretization.  Colors and line styles depict the same information as in previous plots.  Notice that the truncation errors in the numerical methods generated in this work are similar to those of the second order centered difference, that is $O(N^{-2})$.  This is true until $N$ is large enough when the errors begin to grow as $O(N)$.  This is likely due to floating point errors.}\label{fig:Truncation_Error}
\end{figure}

Figure \ref{fig:Truncation_Error} illustrates these errors under refinement.  Initially, the error in each numerical method decays as $O(N^{-2})$, which agrees with the approximation order of the second order centered difference method used to initialize the optimization procedure that determines each weight set.  It is not surprising that these approximations have truncation error that is as good as the method they are initialized from, since the optimization procedure attempts to minimize the objective that, in part, includes the truncation error.  However, it is unfortunate that this procedure does not improve on the approximation order.  For large $N$, the error then grows as $O(N)$ likely due to floating point errors; however, this has not been confirmed.

Further studies on the impact of the truncation error of the method used to initialize the optimization procedure on the order of the numerical method utilizing the optimized weights are warranted.  However, this study is concerned only with the possibility of generating stable methods, whose truncation errors could later be analyzed on problems with known solutions, as is done here.  Still, while these results demonstrate that the numerical methods have truncation error similar to the second order centered difference (which, again, was used as a guess for the iterative process of optimization), it is possible to introduce a set of constraints to enforce a particular truncation error on the methods.  This set of constraints enforces that the weight set is exact for polynomials of a certain degree, e.g.,
\begin{align}
    \sum\limits_{j=1}^{n}w_{k,j}(x_{k,j}-x_{k})^{m}=\mathcal{L}\left((x-x_{k})^{m}\right)|_{x=x_{k}}\nonumber
\end{align}
for $m=0,1,\ldots,m'$, $m'<n-1$.  Since $w_{k,j}$ and $x_{k,j}-x_{k}$ are the same for each value of $k$, enforcement of these constraints leads to only $m'+1$ equations.  The choice of $m'$ creates a numerical method with truncation error of {\color{red}{$O(h_{x}^{m'-1})$}}.  Fixing $m'=n-1$ would use all available decision variables.  This approach of enforcing a particular truncation error on the numerical method was implemented and tested by the author; however, the cases where a stable numerical method were found under these constraints were even fewer than those without the constraints.

\section{Other Significant Observations} \label{sec:otherobs}

In the absence of a concise visual representation of all of the data collected, further observations on the performance of this approach for generating weight sets for approximating the spatial derivative operators when solving the advection-diffusion equation are collected here.
\begin{itemize}
\item When $p=8$, the weight sets generated for $n=3$ deviate little from those in $D'$.  The weights in $D'$ are constructed to have low truncation error at each step, specifically in the case of $n=3$.  Attempts to decrease $J_{\tau q}$, and thus $J$, in this case are more difficult as the initial guess for the quasi-Newton method, which is the weight set used in $D'$, is already likely at or near a local minimum of $J$.  It may be possible to determine weights that both deviate from those in $D'$ and improve upon the value of $J$; however, there is no clear systematic way to determine initial guesses that would do so.
\item In line with the previous observation, generation of weight sets that are both more accurate and stable appears more likely as $n$ increases.  This should be expected as the search space increases in size.
\item Increases in $h_{t}$ beyond $h_{t,N,\nu,s}$ inflate the number of scaled eigenvalues of $D'$ that fall outside the stability region while also increasing their modulus.  Therefore, in practice the number of weight sets that are found that lead to stable numerical methods decreases with the increases in $h_{t}$ examined here.  This is unfortunate as it would be desirable to determine weight sets for which stability is possible even for large values of $h_{t}$ to improve the efficiency of numerical simulation.
\item While there is a maximum value of $h_{t}$ for which the centered difference method is stable for each choice $\nu$ when $s=3$ (but not for $s=2$), this approach to constructing weight sets that lead to stable numerical methods is more successful when $s=2$.  Figure \ref{fig:StabRegionEvals} illustrates that the area of the stability region of AB-3 is larger than that of AB-2, which may account for the greater success.  However, this has not been systematically verified.  Further, AB-3 is a more accurate method than AB-2 and the discussions on accuracy above should analogously apply here.
\end{itemize}

\section{Conclusions} \label{sec:conclusions}

The computational results presented in sections \ref{sec:keyobs} and \ref{sec:otherobs} demonstrate that stable solutions of the advection-diffusion equation can be found, albeit in a limited setting, through the minimization of a recurrent loss function for parameter sets where traditional numerical methods may fail. However, the observations of the same sections highlight an important difficulty with data-driven approaches to solving PDEs--the results can be unpredictable.  The PDE and network model used here can be fully characterized, with known exact solutions to the PDE and with weights in the network that can be easily interpreted as entries in a matrix that represents a discrete approximation of a linear operator.  Yet, even in this simplified case, the results do not suggest the existence of a reliable and robust set of criteria for constructing weight sets that lead to stable numerical methods for solving \eqref{eq:ContinousProblem}.  Most importantly, it does not appear that increasing any of
\begin{itemize}
\item the number of time steps included in the recurrent loss,
\item the smoothness of the initial data,
\item the amount of training data, or
\item the number of iterations in the optimizations method
\end{itemize}
has a predictable impact on determining a desirable set of weights.  

Still, the successes demonstrated here are worth further investigation, and there are many opportunities for further research into stable numerical methods for time-dependent PDEs that rely on the techniques of machine learning.  Some of these include
\begin{itemize}
\item approaches for the selection of initial weight vectors for the quasi-Newton method
\item consideration of entirely different methods of optimization like steepest descent (see, e.g., \cite{Nocedal_Wright_2006})
\item characterization of the size of, e.g., ${\boldsymbol\delta}_{\tau q}^{C}(lh_{t};\boldsymbol{\omega})$, as $Q$ increases and for different value of $l$ to identify potential vanishing gradients
\item alternative network architectures, like long short-term memory (LSTM) \cite{LSTM}, in the case of vanishing gradients
\item allowance for the multi-step weights $\alpha_{l}$ to also vary in the process of minimizing $J$.
\end{itemize}

This approach can also be generalized to other types of discretizations of the spatial and time(-like) variables and to different PDEs.  However, when considering the generation of weights outside of the context of finite differences on uniformly spaced data, their will likely be many more parameters to optimize leading to a more computationally intensive computational procedure to determine the weights.  This is because, for instance, in the presence of data that is not uniformly spaced (e.g., scattered data) a different set of weights will need to be generated for approximating the action of $\mathcal{L}$ at each node.  Similarly, if the problem being solved does not have a periodic solution, then enforcement of boundary conditions is likely to introduce the need for a larger parameter set that needs to be optimized.  In any case, the network representation and backpropagation approach discussed in sections \ref{sec:networkrep} and \ref{sec:backprop} would need to be altered to accommodate these new parameters and the way that they appear in the problem.

\section{Funding}

This work was supported in part by the Air Force Office of Scientific Research under the project ``Kernel Methods with Machine Learning and Adaptivity".  The author has no other relevant financial or non-financial interests to disclose.

\bibliographystyle{unsrt}      
\bibliography{references}   

\end{document}